\pdfoutput=1
\RequirePackage{ifpdf}
\ifpdf 
\documentclass[pdftex]{sigma}
\else
\documentclass{sigma}
\fi

\usepackage{mathtools}
\usepackage[all,cmtip]{xy}

\usepackage{tikz, pgfplots}

\usepackage{floatrow}
\usetikzlibrary{graphs,quotes,decorations.pathmorphing, decorations.markings}

\DeclareMathOperator{\tr}{tr}

\DeclareMathOperator{\ind}{ind}

\DeclareMathOperator{\id}{id}

\DeclareMathOperator{\diag}{diag}

\DeclareMathOperator{\ran}{ran}

\DeclareMathOperator{\ad}{ad}
\DeclareMathOperator{\sgn}{sgn}

\newcommand{\eq}{\mathbin{\rotatebox[origin=c]{90}{$=$}}}

\newtheorem{Theorem}{Theorem}[section]
\newtheorem{Corollary}[Theorem]{Corollary}
\newtheorem{Lemma}[Theorem]{Lemma}
\newtheorem{Proposition}[Theorem]{Proposition}
 { \theoremstyle{definition}
\newtheorem{Definition}[Theorem]{Definition}

\newtheorem{Example}[Theorem]{Example}
\newtheorem{Remark}[Theorem]{Remark} }

\numberwithin{equation}{section}

\newtheorem*{Theorem*}{Theorem}

\begin{document}

\allowdisplaybreaks

\newcommand{\arXivNumber}{2307.13324}

\renewcommand{\PaperNumber}{022}

\FirstPageHeading

\ShortArticleName{Boundary Value Problems for Dirac Operators on Graphs}

\ArticleName{Boundary Value Problems for Dirac Operators\\ on Graphs}

\Author{Alberto RICHTSFELD}

\AuthorNameForHeading{A.~Richtsfeld}

\Address{Institut f\"ur Mathematik, Universit\"at Potsdam, D-14476, Potsdam, Germany}
\Email{\href{mailto:richtsfeld@uni-potsdam.de}{richtsfeld@uni-potsdam.de}}

\ArticleDates{Received August 21, 2023, in final form February 28, 2024; Published online March 19, 2024}

\Abstract{We carry the index theory for manifolds with boundary of B\"ar and Ballmann over to first order differential operators on metric graphs. This approach results in a short proof for the index of such operators. Then the self-adjoint extensions and the spectrum of the Dirac operator on the complex line bundle are studied. We also introduce two types of boundary conditions for the Dirac operator, whose spectrum encodes information of the underlying topology of the graph.}


\Classification{34B45; 58J50}

\section{Introduction}

Schr\"odinger operators on metric graphs, known as quantum graphs, have a long history of study with applications in many fields, including mathematics, physics and engineering. They provide toy models for quantum mechanics, but also describe the motion of particles through thin metallic wires called nanowires. We will mainly focus on the following mathematical aspects of the theory of quantum graphs, namely self-adjoint extensions, which are realised by imposing boundary conditions at each vertex, and the index theory of these operators. A full treatment of self-adjoint extensions can be found in \cite{FullingKuchmentWilson, Harmer, KostrykinSchrader, Kuchment2014}, which is summarised by Kuchment~\cite{kuchment2008quantum}, the take-away being that for the Schr\"odinger operator self-adjoint extensions are given by a~combination of Dirichlet, von Neumann and Robin boundary conditions. Index theory for quantum graphs was treated by Fulling, Kuchment and Wilson \cite{FullingKuchmentWilson}, relating the index to heat kernel asymptotics. A nice overview of the study of quantum graphs is given for example by Kuchment~\cite{kuchment2008quantum}, the monograph \cite{BerKuch} by Berkolaiko and Kuchment gives an introduction to the theory of quantum graphs.

In order to study the effect of spin on the quantum mechanics of graphs, the Dirac operator has to be considered.
The first ones to study a two-dimensional, quantum-mechanical Dirac operator on graphs were Bulla and Trenkler \cite{bulla}. They studied self-adjoint extensions of the Dirac operator on directed graphs and computed its spectrum on a~junction of three wires, a graph for which the Schr\"odinger operator had previously been considered by Exner and \v{S}eba~\cite{exner}. Bolte and Harrison \cite{bolte} considered this Dirac operator in a similar way to the approach taken in this paper: the graph under consideration is metric in the sense that each edge has an associated length, and its boundary conditions are given by a unitary matrix describing complex transition probability amplitudes. Post \cite{post} introduced a discrete notion of the Dirac operator and computed the index of the discrete and metric Dirac operators, by finding that it is equal to the Euler-characteristic of the graph. For more information on the study of Dirac operators on metric graphs, see the survey article \cite{harrisonsurv} written by Harrison.

The index theory of Dirac operators on manifolds with boundary was initiated by the seminal paper of Atiyah, Patodi and Singer \cite{aps}, whose main achievement was the discovery of a non-local boundary condition, formulated in terms of the spectrum of the boundary operator.
The theory centred on boundary value problems for Dirac operators was then simplified and unified by B\"ar and Ballmann \cite{bb, bb2}, in which they characterised all extensions of the Dirac operators with linear subspaces of the Czech-space, a Banach space formed from the eigenspaces of the boundary operator. These results were then extended by B\"ar and Bandara \cite{lashi} to general first-order differential operators. We will work in the spirit of the B\"ar--Ballmann framework, but the theory does not apply directly to the one-dimensional case, since there is no notion of a~boundary operator in this case, since the boundary is zero-dimensional.

The point of view in this paper is as follows: A metric graph can be reinterpreted as a~one-dimensional Riemannian manifold with boundary, where each connected component is isomorphic to a closed interval. We then consider the Dirac operator on this one-dimensional manifold, subject to boundary conditions that respect the underlying graph structure. The aim is to study the self-adjoint extensions and the index of this operator, and to find interesting boundary conditions.

After fixing notation, we consider general first-order differential operators on metric graphs with coefficients in a general bundle, without specifying the boundary conditions. The B\"ar--Ballmann approach leads to a new proof of the index theorem for first-order operators on graphs (see Section~\ref{section:index} and Theorem~\ref{thm:ind}), which has already been proved by different methods for Dirac operators by Post \cite{post} and for general operators by Fulling, Kuchment and Wilson \cite{FullingKuchmentWilson}.

We then turn to the scalar Dirac operator. We first study its self-adjoint extensions and find that these are given by boundary conditions which are graphs of unitary endomorphisms, see Theorem~\ref{thm:unit}. Self-adjoint extensions have already been studied by Bulla and Trenkler \cite{bulla} and Bolte and Harrison \cite{bolte} for another Dirac operator acting on a two-dimensional bundle over the graph. Since the scalar Dirac operator acts on a one-dimensional bundle, the conditions on self-adjoint extensions are more restrictive, leading to fewer but simpler such extensions. In particular, self-adjoint extensions exist only if the graph is a directed Eulerian graph.

In Section~\ref{sec:spectrum}, we take a closer look at the spectrum of the scalar Dirac operator under boundary conditions given by graphs of endomorphisms. Here, the endomorphisms are not necessarily unitary, so that the spectrum is no longer real. Theorem~\ref{thm:spec} tells us that the zero-locus of the multidimensional characteristic polynomial of the endomorphism whose graph gives the boundary condition completely determines the spectrum of the Dirac operator. In the case that all lengths are commensurable, the spectrum can be computed explicitly in terms of the eigenvalues of the endomorphism. For a reducible boundary condition, we obtain a splitting property of the spectrum, i.e., we can find subgraphs such that the spectrum of the scalar Dirac operator is the union of the spectra of the scalar Dirac operator of these subgraphs (see Theorem~\ref{prop:decomp}).

In the last part, we discuss two types of boundary conditions. The first is given by certain permutation matrices which encode decompositions of the underlying graph into directed trails, the spectrum of the Dirac operator is then explicitly calculable and is fully determined by the lengths of the trails appearing in the decomposition. The second boundary condition is given by the directed adjacency matrix, which distributes the inflow at a vertex evenly among the outgoing edges. This boundary condition is in general no longer self-adjoint, but it links the spectrum of the Dirac operator to the appearance of cycles in the graph. Thus, both types of boundary conditions reveal properties of the cycles and trails present in the graph, linking its topology to the spectrum of the Dirac operator.

\section{Preliminaries}
A metric graph $G=(V,E,l)$ is given by a directed multigraph with vertex set $V$ and edge set~$E$ and a length function $l\colon E\rightarrow (0,\infty)$, $e\mapsto l_e$. For each edge $e\in E$ we denote by $\partial_-e$ the tail of~$e$ and by $\partial_+e$ the head of $e$. For our purposes, we will always assume that our metric graph is finite and has no isolated vertices.

We associate to $G$ the following disconnected (Riemannian) $1$-manifold with boundary
\[
	\hat{G}=\bigsqcup_{e\in E}[0,l_e],
\]
on which we have the following function spaces:
\begin{gather*}
	L^2\bigl(\hat{G},\mathbb{C}^r\bigr)= \bigoplus_{e\in E} L^2([0,l_e], \mathbb{C}^r),\\
 H^1\bigl(\hat{G},\mathbb{C}^r\bigr)=\bigoplus_{e\in E}H^1([0,l_e],\mathbb{C}^r),\\
	C^\infty\bigl(\hat{G},\mathbb{C}^r\bigr)=\bigoplus_{e\in E} C^\infty ([0,l_e],\mathbb{C}^r).
\end{gather*}
Next, we define a vector space in which we store boundary values of functions. As each function has two boundary values for each edge, one at its tail and one at its head, we need a vector space that is generated by two copies of the edge set, denoted by $\bar{E}= E_+\sqcup E_-$, where
$E_\pm := \{e_\pm\mid e\in E\}$.
Define the said vector space to be
$\mathbb{C}^{r,\bar{E}}= \mathbb{C}^r\otimes \mathbb{C}^{\bar{E}}$.
Standard Sobolev estimates tell us that each function in $H^1\big(\hat{G},\mathbb{C}^r\big)$ is continuous when restricted to a single edge and that there is a well-defined bounded operator
\[
	\tr \colon \ H^1\big(\hat{G},\mathbb{C}^r\big)\rightarrow \mathbb{C}^{r,\bar{E}},\qquad f=(f_e)_{e\in E}\mapsto \sum_{e\in E}f_e(0)\otimes e_-+f_e(l_e)\otimes e_+.
\]
This allows us to introduce boundary conditions on the graph: A boundary condition is a~linear subspace $B\subseteq \mathbb{C}^{r,\bar{E}}$, we then denote the space of $H^1$-functions which fulfill the boundary condition $B$ by
\[
	H^1_B\bigl(\hat{G},\mathbb{C}^r\bigr)=\bigl\{f\in H^1\bigl(\hat{G},\mathbb{C}^r\bigr)\mid \tr(f)\in B\bigr\}.
\]
If $B=\{0\}$, we write $H^1_0\bigl(\hat{G},\mathbb{C}^r\bigr)$ instead.

There exists a boundary map
$\partial \colon \bar{E}\rightarrow V$,
which for an edge $e\in E$ sends $e_-\in\bar{E}$ to the tail
$\partial_-e\in V$ and $e_+\in \bar{E}$ to the head $\partial_+ e\in V$ of~$e$.
Given a vertex $v\in V$, set
$\bar{E}_v := \partial^{-1}(\{v\})\subseteq \bar{E}$,
which in turn generates a subspace $\mathbb{C}^{r,\bar{E}_v}$ of $\mathbb{C}^{r,\bar{E}}$.
\begin{Definition}
A boundary condition $B$ is {\it local} iff there are linear subspaces $B_v\subseteq \mathbb{C}^{r,\bar{E}_v}$, such that
$B=\bigoplus_{v\in V}B_v$.
\end{Definition}

\begin{Remark}
	This notion of locality is different from the notion of locality introduced by B\"ar and Ballmann \cite[Definition 7.9]{bb}. In fact, every boundary condition on the manifold $\bigsqcup_{e\in E}[0,l_e]$ is local in the sense of B\"ar and Ballmann, and non-locality is a phenomenon that occurs only in higher dimensions. Locality in our sense means that the boundary condition sees the structure of the graph and imposes a condition per vertex such that the conditions at two different vertices are independent of each other.
\end{Remark}

We define also
$\mathbb{C}^{r,E}:= \mathbb{C}^r\otimes \mathbb{C}^E$,
as well as the boundary maps $\partial_\pm \colon E\rightarrow V$, mapping $e \rightarrow \partial_\pm e$. The subsets $E_v^\pm := \partial_\pm^{-1}(\{v\})\subseteq E$ induce subspaces
$\mathbb{C}^{r,E_v^\pm}\subseteq \mathbb{C}^{r,E}$.
Note that there are canonical inclusions $\iota_\pm\colon \mathbb{C}^{r,E}\rightarrow \mathbb{C}^{r,\bar{E}}$ sending $e$ to $e_\pm$. For $v\in V$, the subspace $\mathbb{C}^{r,\bar{E}_v}$ is equal to the direct sum $\iota_+\big(\mathbb{C}^{r,E_v^+}\big)\oplus \iota_-\big(\mathbb{C}^{r,E_v^-}\big)$.

In the case that $r=1$, we omit $r$ in the above definitions. We equip all so far defined vector spaces with the canonical inner products making any two elements of the generating sets orthonormal to each other.

\begin{Example}
	Let $G$ be a graph consisting of one edge $e$ of length $l$, compare Figure \ref{figure:spaces}. We then have
	\begin{gather*}
		\mathbb{C}^{E} = \mathbb{C}\cdot e,\qquad
		\mathbb{C}^{\bar{E}} = \mathbb{C}\cdot e_-\oplus \mathbb{C}\cdot e_+,\qquad
		C^\infty\bigl(\hat{G},\mathbb{C}\bigr) = C^\infty([0,l],\mathbb{C}).
	\end{gather*}
	For $f\in C^\infty\bigl(\hat{G},\mathbb{C}\bigr)$,
	\[
		\operatorname{tr}(f)= f(0)\cdot e_- + f(l)\cdot e_+
	\]
	holds.

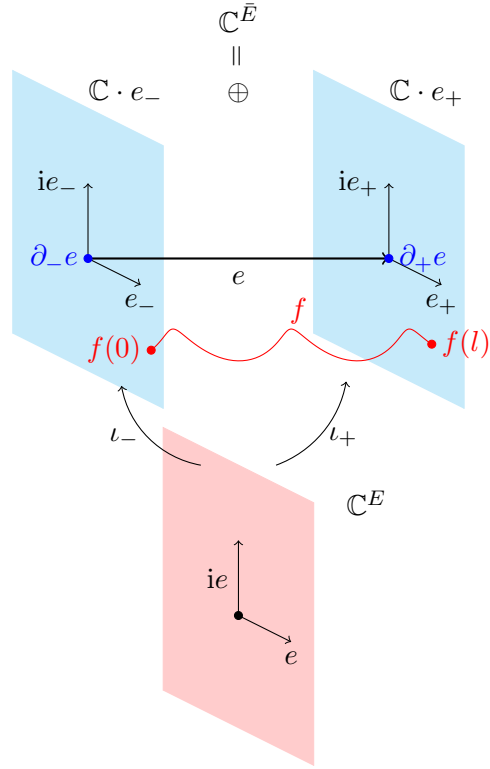
\begin{figure}[t]
	\centering
	\begin{tikzpicture}[domain=0:4, samples=100]
		\node (up) at (3,3.2) {$\mathbb{C}^{\bar{E}}$};
¸		\node (eq) at (3,2.7) {$\eq$};
		\node (plus) at (3,2.2) {$\oplus$};
		
		\filldraw[cyan!20] (0,-1) -- (0,2.5) -- (2,1.5) -- (2,-2);
		\node (N1) at (1.5,2.2) {$\mathbb{C}\cdot e_-$};
		\draw[->] (1,0) -- (1.7,-0.35);
		\node[anchor=north] (e-) at (1.7,-0.35) {$e_-$};
		\draw[->] (1,0) -- (1,1);
		\node[anchor=east] (ie-) at (1,1) {${\rm i}e_-$};
		
		\filldraw[cyan!20] (4,-1) -- (4,2.5) -- (6,1.5) -- (6,-2);
		\node (N2) at (5.5,2.2) {$\mathbb{C}\cdot e_+$};
		\draw[->] (5,0) -- (5,1);
		\node[anchor=east] (ie+) at (5,1) {${\rm i}e_+$};
		\draw[->] (5,0) -- (5.7,-0.35);
		\node[anchor=north] (e+) at (5.7,-0.35) {$e_+$};
		
		\draw[->,thick] (1,0) -- (5,0);
		\filldraw[blue] (1,0) circle (1.5pt) node[anchor=east]{$\partial_- e$};
		\filldraw[blue] (5,0) circle (1.5pt) node[anchor=west]{$\partial_+ e$};
		\node (edge) at (3,-0.25) {$e$};
		
		\draw[color=red] plot ({\x+1+0.7*(0.2*cos(4*\x r)+1)},{(0.2*sin(4*\x r)-0.8)-0.35*(0.2*cos(4*\x r)+1)});
		\filldraw[red] (1.84,-1.22) circle (1.5pt) node[anchor=east]{$f(0)$};
		\filldraw[red] (5.57,-1.14) circle (1.5pt) node[anchor=west]{$f(l)$};
		\node[anchor=south, red] (f) at (3.8,-1) {$f$};
		
		\node (N3) at (4.7,-3.25) {$\mathbb{C}^{E}$};
		\filldraw[red!20] (2,-5.75) -- (2,-2.25) -- (4,-3.25) -- (4,-6.75);
		\filldraw[black] (3,-4.75) circle (1.5pt);
		\draw[->] (3,-4.75) -- (3.7,-5.1);
		\node[anchor=north] (e) at (3.7,-5.1) {$e$};
		\draw[->] (3,-4.75) -- (3,-3.75);
		\node[anchor=east] (ie) at (3,-4.25) {${\rm i}e$};
		
		\draw[->] (2.5,-2.75)arc(260:190:1.3) node[pos=0.5, anchor=east] {$\iota_-$};
		\draw[->] (3.5,-2.75)arc(290:350:1.45) node[pos=0.45, anchor=west] { $\iota_+$};
	\end{tikzpicture}
	\caption{The spaces $\mathbb{C}^{\bar{E}}$, $\mathbb{C}^E$ for a graph consisting of a single edge $e$. The continuous function $f$ has boundary values in $\mathbb{C}^{\bar{E}}$.}
	\label{figure:spaces}
\end{figure}
	
	There are two cases: Either we have $\partial_-e \neq \partial_+ e$, or $\partial_-e= \partial_+ e$, in which case $G$ is a loop.
	In the first case, we have
	\begin{alignat*}{3}
		&\mathbb{C}^{\bar{E}_{\partial_- e}}= \mathbb{C}\cdot e_-,\qquad && \mathbb{C}^{\bar{E}_{\partial_+ e}}= \mathbb{C}\cdot e_+,&\\
		&\mathbb{C}^{E^-_{\partial_- e}} = \mathbb{C}\cdot e, && \mathbb{C}^{E^-_{\partial_+ e}}= \{0\},&\\
		&\mathbb{C}^{E^+_{\partial_- e}} = \{0\}, && \mathbb{C}^{E^+_{\partial_+ e}} = \mathbb{C}\cdot e,&
	\end{alignat*}
	and the local boundary conditions are given by the subspaces
	\[
		\{0\},\  \mathbb{C}^{\bar{E}_{\partial_- e}},\  \mathbb{C}^{\bar{E}_{\partial_+ e}},\  \mathbb{C}^{\bar{E}}\subseteq \mathbb{C}^{\bar{E}}.
	\]
	In the second case $\partial_-e= \partial_+ e=v$, we have
	\begin{alignat*}{4}
		&\mathbb{C}^{\bar{E}_v}=\mathbb{C}^{\bar{E}},\qquad &&
		\mathbb{C}^{E^-_v}=\mathbb{C}^{E},\qquad &&
		\mathbb{C}^{E^+_v}=\mathbb{C}^{E}&
	\end{alignat*}
	and the local boundary conditions are given by any linear subspace $B$ of $\mathbb{C}^{\bar{E}}$.
\end{Example}

Given two endomorphism fields $\sigma, \mathbb{V}\in C^\infty\big(\hat{G},\mathbb{C}^{r\times r}\big)$, we can define a first order differential operator
\[
	D\colon \ C^\infty\bigl(\hat{G},\mathbb{C}^r\bigr)\rightarrow C^\infty \bigl(\hat{G},\mathbb{C}^r\bigr), \qquad f\mapsto \sigma f' +\mathbb{V}f.
\]
We will always assume that the differential operator $D$ is elliptic, i.e., that $\sigma$ is non-singular at every point. $D$ naturally extends to a bounded operator
\[
	D_{\max}\colon \ H^1\bigl(\hat{G},\mathbb{C}^r\bigr)\rightarrow L^2\bigl(\hat{G},\mathbb{C}^r\bigr)
\]
 and consequently for a boundary condition $B$ we can restrict $D_{\max}$ to $H^1_B\bigl(\hat{G},\mathbb{C}^r\bigr)$, yielding the operator $D_B$. If $B=\{0\}$, write $D_{\min}$ instead. These operators will be the center of our study.

It is straightforward to check that the unbounded operator
\[
	D_B\colon \ H^1_B\bigl(\hat{G},\mathbb{C}^r\bigr)\subseteq L^2\bigl(\hat{G},\mathbb{C}^r\bigr)\rightarrow L^2\bigl(\hat{G},\mathbb{C}^r\bigr)
\]
is closed, densely defined and a Fredholm operator. Furthermore, we have the following Green's formula.
\begin{Lemma}
	For $f,g\in H^1\bigl(\hat{G},\mathbb{C}^r\bigr)$, we have that
	\[
		\langle D_{\max}f,g\rangle - \bigl\langle f, D_{\max}^\dagger g\bigr\rangle = \langle \sigma_0 \tr(f),\tr(g)\rangle,
	\]
	where $D^\dagger$ is the formal adjoint of $D$ and
	\[
		\sigma_0\colon \ \mathbb{C}^{r,\bar{E}}\rightarrow \mathbb{C}^{r,\bar{E}}, \qquad \sum_{e\in E} a_e\otimes e_- + b_e\otimes e_+\mapsto \sum_{e\in E} -\sigma_e(0)a_e\otimes e_- + \sigma_e(l_e)b_e \otimes e_+,
	\]
	where $\sigma = (\sigma_e)_{e\in E}$.
\end{Lemma}

We directly obtain the following assertion.

\begin{Proposition}\label{prop:sa}
	Let $B$ be a boundary condition, define the adjoint boundary condition by
	\[
		B^{\mathrm{ad}}:=\bigl\{y\in \mathbb{C}^{r,\bar{E}}\mid \forall b\in B\colon \langle \sigma_0 b, y\rangle =0\bigr\}=(\sigma_0 B)^\perp.
	\]
	Then, the adjoint of $D_B$ is given by $D_{B^{\mathrm{ad}}}^\dagger$. In particular, $D_B$ is self-adjoint iff $D$ is formally self-adjoint and $B=B^{\mathrm{ad}}$.
\end{Proposition}
\begin{Corollary}
	The adjoint of $D_{\max}$ is $D^\dagger_{\min}$ and the adjoint of $D_{\min}$ is $D^\dagger_{\max}$.
\end{Corollary}

\section[The index of D\_B]{The index of $\boldsymbol{D_B}$}\label{section:index}

By elementary ODE-theory we can solve the equation
\[
	D_{\max}f=g
\]
on each edge separately for any $g\in L^2\bigl(\hat{G},\mathbb{C}^r\bigr)$. Hence $D_{\max}$ is surjective. Furthermore, it is easily seen that the dimension of the solution space of the equation
\[
	D_{\max}f=0
\]
is equal to $r$ on each single edge, and since the kernel of $D_{\max}$ is just the direct sum of the kernels of the single edges, we obtain
\[
	\dim\ker D_{\max} = r\lvert E\rvert.
\]
We have therefore proven the lemma.

\begin{Lemma}
	The following equations hold:
	\[
		\ind D_{\max} = r\lvert E\rvert, \qquad \ind D_{\min}= - r\lvert E\rvert.
	\]
\end{Lemma}

In order to prove the general index formula, we pursue the approach developed by B\"ar and Ballmann. We need the following proposition:

\begin{Proposition}\label{prop:bb}
	Let $H$ be a Hilbert space, $E$ and $F$ Banach spaces and let $L\colon H\rightarrow E$ and $P\colon H\rightarrow F$ be bounded linear maps. Assume that $P:H\rightarrow F$ is onto. Then, $L|_{\ker P}\colon \ker P\rightarrow E$ is Fredholm of index $k$ if and only if $L\oplus P\colon H\rightarrow E\oplus F$ is Fredholm of index $k$.
\end{Proposition}
\begin{proof}
	See \cite[Proposition A.1]{bb}.
\end{proof}

Analogous to \cite[Corollary 8.7]{bb}, we have the following.

\begin{Lemma}\label{cor:fred}
	Let $B$ be a subspace of $\mathbb{C}^{r,\bar{E}}$, $B^\perp$ be the orthogonal complement to $B$ and let $P\colon \mathbb{C}^{r,\bar{E}}\rightarrow\mathbb{C}^{r,\bar{E}}$ be the orthogonal projection onto $B^\perp$. Then,
	\begin{align*}
		L\colon \ H^1\bigl(\hat{G},\mathbb{C}^r\bigr)&\rightarrow L^2\bigl(\hat{G},\mathbb{C}^r\bigr)\oplus B^\perp\\
		g&\mapsto (D_{\max}g, P\tr g)
	\end{align*}
	is a Fredholm operator with the same index as $D_B$.
\end{Lemma}
\begin{proof}
	The kernel of $P\tr$ is equal to $H^1_B\bigl(\hat{G},\mathbb{C}^r\bigr)$. The statement then follows directly from Proposition \ref{prop:bb}.
\end{proof}

As in \cite[Corollary 8.8]{bb}, we go on to obtain the lemma.

\begin{Lemma}\label{lem:deform}
	Let $B_1\subseteq B_2 \subseteq \mathbb{C}^{r,\bar{E}}$ be boundary conditions. Then,
	\[
		\ind D_{B_2}=\ind D_{B_1}+ \dim B_2/B_1.
	\]
\end{Lemma}
\begin{proof}
	Since $B_1$ is a subspace of $B_2$,
	the orthogonal complement $B_2^\perp$ is a subspace of $B_1^\perp$. Denote the inclusion of $B_2^\perp$ into $B_1^\perp$ by $\iota$. The diagram
	\[
		\xymatrix{
			& L^2\bigl(\hat{G},\mathbb{C}^r\bigr)\oplus B_2^\perp \ar@{^{(}->}[dd]^{\id \oplus\, \iota} \\
			H^1\bigl(\hat{G},\mathbb{C}^r\bigr) \ar[ru]^{(D,P_2)} \ar[rd]^{(D,P_1)} \\
			& L^2\bigl(\hat{G},\mathbb{C}^r\bigr)\oplus B_1^\perp
		}
	\]
	commutes, where, as in Lemma \ref{cor:fred}, the operator $P_i$ is the composition of $\tr$ with the orthogonal projection onto $B_i^\perp$. Clearly, $\id\oplus\,\iota$ is a Fredholm operator with index
	\[
		\ind (\id\oplus\,\iota)=\ind (\iota) = -\dim\bigl(B^\perp_1/B^\perp_2\bigr) = -\dim(B_2/B_1).
	\]
	Since the index is additive, we have
	\begin{align*}
		\ind D_{B_1}= \ind (D,P_1)
		= \ind (D,P_2) + \ind \id \oplus \iota
		= \ind D_{B_2} -\dim B_2/B_1.\tag*{\qed}
\end{align*}
\renewcommand{\qed}{}
\end{proof}

\begin{Theorem}\label{thm:ind}
	For any boundary condition $B$, the following formula for the index of $D_B$ holds:%
	\begin{equation*}
		\ind(D_B)=\dim B -r|E|.
	\end{equation*}
\end{Theorem}

\begin{proof}
	Applying Lemma \ref{lem:deform}, the index of $D_B$ is given by
	\begin{align*}
		\dim D_B = \dim B + \ind D_0
		= \dim B + \ind D_{\min}
		= \dim B -r|E|.\tag*{\qed}
\end{align*}\renewcommand{\qed}{}
\end{proof}

\section{The scalar Dirac operator}

For the rest of the paper, we will focus on the scalar Dirac operator, that is, we take $r=1$ and set $D= {\rm i}\frac{{\rm d}}{{\rm d}x}$.
We will first study its self-adjoint extensions.

\subsection{Self-adjoint extensions}

Let $A\colon \mathbb{C}^E\rightarrow \mathbb{C}^E$ be a linear map, we can then define an associated boundary condition
\[
\Gamma(A):=\bigl\{\iota_+(x)+\iota_-(Ax)\mid x\in\mathbb{C}^E\bigr\}\subseteq \mathbb{C}^{\bar{E}}
\]
which we call the {\it graph of $A$}. We say that $A$ is a {\it $G$-endomorphism}, if $\Gamma(A)$ is a local boundary condition. In this case $A$ restricts to a homomorphism
\[
	A_v \colon \ \mathbb{C}^{E_v^+}\rightarrow \mathbb{C}^{E_v^-}
\]
for all $v\in V$.

\begin{Theorem}\label{thm:unit}
	Let $B$ be a boundary condition on $G$. Then, $D_B$ is self-adjoint iff there exists a~unitary map $A\colon \mathbb{C}^E\rightarrow\mathbb{C}^E$, such that $B=\Gamma(A)$.
\end{Theorem}
\begin{proof}
	Let $B$ be such that $D_B$ is self-adjoint. Since the index of a self-adjoint operator is zero, by Theorem~\ref{thm:ind}, the dimension of $B$ is equal to the number of edges, denoted by $n$. Choosing a basis $b_1,\dots,b_n$ of $B$, let
	\[
		b_i^\pm=\pi_\pm(b_i)\in \mathbb{C}^E
	\]
	for $i=1,\dots,n$, where $\pi_\pm\colon \mathbb{C}^{\bar{E}}\rightarrow \mathbb{C}^E$ is the map sending $e_\pm$ to $e$ and $e_\mp$ to $0$. Now, assume that there exist $\lambda_1,\dots,\lambda_n$, such that
	\[
		\sum_{i=1}^n \lambda_i b^+_i =0.
	\]
	Letting
	\[
		b=\sum_{i=1}^n \lambda_i b_i,
	\]
	we arrive at $\pi_+(b)=0$. Since $B$ is equal to its adjoint boundary condition, by Proposition \ref{prop:sa} we have
	\begin{align*}
		0= \langle \sigma_0(b),b\rangle={\rm i} \langle \pi_{-}(b),\pi_{-}(b)\rangle -{\rm i}\langle \pi_{+}(b),\pi_{+}(b)\rangle = {\rm i} \langle \pi_{-}(b),\pi_{-}(b)\rangle.
	\end{align*}
	Therefore, $\pi_-(b)$ is also zero, which implies that $b$ is equal to zero. Since $b_1,\dots,b_n$ is a basis of~$B$,
	\[
		\lambda_1=\dots =\lambda_n=0
	\]
	follows. Consequently, the vectors $b^+_1,\dots,b^+_n$ are linearly independent.
	Considering that the dimension of $\mathbb{C}^E$ is equal to $n$, it follows that $b^+_1,\dots,b^+_n$ forms a basis of $\mathbb{C}^E$.
Define~$A$ as the linear map sending $b^+_i$ to $b^-_i$ for $i\in\{1,\dots,n\}$. Since $B=B^{\ad}$, we have
	\[
\bigl\langle \sigma_0(b_i),b_j\bigr\rangle=0
	\]
	for $i,j\in\{1,\dots,n\}$, hence
	\[
{\rm i} \bigl\langle b^-_i,b^-_j\bigr\rangle -{\rm i}\bigl\langle b^+_i,b^+_j\bigr\rangle=0.
	\]
	Rearranging this equation and using the definition of $A$, this yields
	\begin{align*}
		\bigl\langle b^+_i,b^+_j\bigr\rangle=\bigl\langle Ab^+_i,Ab^+_j\bigr\rangle,
	\end{align*}
	which proves that $A$ is unitary.
	It is easily checked that $A$ is a $G$-endomorphism. Therefore,
	each $b_i$ is contained in $\Gamma(A)$ for $i\in\{1,\dots,n\}$ and, consequently, the same can be said for $B$. Equality, then, holds by equality of dimensions.
	
	Conversely, it is easily checked that for a given unitary map $A\colon \mathbb{C}^E\rightarrow \mathbb{C}^E$, the graph $\Gamma(A)$ is a boundary condition such that $\Gamma(A)^{\ad}=\Gamma(A)$.
\end{proof}

Unfortunately, the existence of local, self-adjoint boundary conditions poses quite some restrictions on the graph itself:

\begin{Proposition}\label{prop: s.a. implies eulerian}
	If there exists a local boundary condition $B$ such that $D_B$ is self-adjoint, then every connected component of $G$ is an Eulerian directed multigraph.
\end{Proposition}

\begin{proof}
	By Theorem~\ref{thm:unit}, there exists a unitary $G$-endomorphism $A$ such that $B$ is the graph of $A$. Restricting $A$ to each vertex space, we then obtain isometries $A_v$ from
	$\mathbb{C}^{E^+_v}$ to $\mathbb{C}^{E^-_v}$, and hence the dimensions of the two spaces must be equal, implying that each vertex has as many ingoing as outgoing edges.
\end{proof}

\subsection{The spectrum of the Dirac operator}\label{sec:spectrum}

We proceed with the study of the spectrum $\sigma(D_B)$ of $D_B$ and start with some general properties:

\begin{Proposition}\label{prop:ptspctrm}
	Let $B$ be a boundary condition on $G$. If the dimension of $B$ is not equal to the number of edges $n$ of $G$, then the spectrum of $D_B$ is equal to $\mathbb{C}$. Otherwise, the spectrum of~$D_B$ is equal to the point spectrum, i.e., $\sigma(D_B)=\sigma_p(D_B)$.
\end{Proposition}

\begin{proof}
	Assume that $\lambda$ is an element of the resolvent set of $D_B$. Then the operator $D_B-\lambda I$ is bijective, and hence the index of $(D-\lambda I)_B= D_B -\lambda I$ must be zero. Theorem~\ref{thm:ind} implies that the dimension of $B$ is equal to $n$. Consequently, if the dimension of $B$ is not equal to $n$, the resolvent set of $D_B$ is empty.
	
	Now, let the dimension of $B$ be equal to $n$. For $\lambda\in\mathbb{C}$, we need to show that injectivity of $D_B-\lambda I$ implies surjectivity. Assuming that $D_B-\lambda I$ is injective,
	by Theorem~\ref{thm:ind}, the codimension of $ \ran (D_B -\lambda I) $ is zero. Furthermore, $D_B-\lambda I$ has closed range since it is Fredholm. Together, this implies that $D_B-\lambda I$ is surjective.
\end{proof}

We now turn to boundary conditions given by $G$-endomorphisms. The spectrum of the Dirac operator under such boundary conditions is closely related to the eigenvalues of the defining $G$-endomorphism. However, the different lengths of the edges distort the powers of the variables in the characteristic polynomial, such that it is necessary to introduce a multi-dimensional version. To this end, introduce the coordinate functions
\begin{align*}
	x_e \colon \ \mathbb{C}^E\rightarrow \mathbb{C},\qquad x_e\biggl(\sum_{e^\prime\in E}a_{e^\prime}e^\prime \biggr)= a_e.
\end{align*}

\begin{Definition}
	For a $G$-endomorphism $A$, define {\it the multi-dimensional characteristic polynomial $P_A$}
of $A$ to be the polynomial on $\mathbb{C}^E$ given by
	\[
		P_A(x)=\det (\diag(x)-A ),
	\]
	where $\diag(x)$ is the linear map on $\mathbb{C}^E$ given by $\sum_{e\in E} x_e e\otimes {\rm e}^\ast$.
	The \emph{characteristic function} $P_{A,l}\colon \mathbb{C}\rightarrow\mathbb{C}$ with respect to the length function $l$ is given by
	\[
		P_{A,l}(\lambda)=P_A(\exp({\rm i}\lambda l)),
	\]
	where $\exp({\rm i}\lambda l) =\sum_{e\in E} \exp({\rm i}\lambda l_e)e$.
\end{Definition}

The next theorem relates the eigenvalues of the Dirac operator to the zeroes of the characteristic function of $A$.

\begin{Theorem}\label{thm:spec}
	Let $A$ be a $G$-endomorphism. Then the spectrum of $D_{\Gamma(A)}$ is equal to the set of zeroes of $P_{A,l}$, i.e.,
	\begin{equation*}
		\sigma\big(D_{\Gamma(A)}\big)=P_{A,l}^{-1}(\{0\}).
	\end{equation*}
	The multiplicity of an eigenvalue $\lambda$ is given by
	\begin{equation*}
		m_{\Gamma(A)}(\lambda)=\dim\ker(\diag(\exp({\rm i}\lambda l))-A).
	\end{equation*}
\end{Theorem}

\begin{Remark}
	Note that $P_{A,l}(\lambda)$ does not capture the kernel of $A$. For instance, if $A=0$, then the spectrum of $D_{\Gamma(0)}$ is empty.
\end{Remark}

\begin{proof}
	Let $\lambda\in \sigma\bigl(D_{\Gamma(A)}\bigr)$. Since the dimension of $\Gamma(A)$ is equal to the number of edges of $G$, by Proposition \ref{prop:ptspctrm}, there exists an eigenvector $\phi\in C^\infty\big(\hat{G}\big)$ to the eigenvalue $\lambda$, i.e., $\phi$ fulfills
	\[D\phi =\lambda \phi, \qquad \tr \phi \in \Gamma(A).\]
	The first condition implies that on each edge $\phi$ is of the form $w_e \exp(-{\rm i}\lambda x_e)$ for some $w_e\in \mathbb{C}$ and together with the second condition we obtain the equation
	\[A\sum w_e \exp(-{\rm i}\lambda l_e)e= \sum w_e e,\]
which implies $P_{A,l}(\lambda) = 0$. Conversely, it is easy to see
that any element of $\ker(\diag(\exp({\rm i}\lambda l))-A)$ defines an
eigenfunction of $D_{\Gamma(A)}$ to the eigenvalue $\lambda$.
\end{proof}

If all edges of the directed multigraph have the same length, the spectrum of $D_{\Gamma(A)}$ is particularly easy to compute:

\begin{Corollary}
	Let all edges have the same length $l$ and $A$ be a $G$-endomorphism. Let $\mu_1,\dots,\mu_m$ be the non-zero eigenvalues of $A$. Then, there
	exist unique $\alpha_1,\dots,\alpha_m\in\mathbb{R}$ and $\varphi_1,\dots,\varphi_m\in[0,2\pi)$ such that
	\[\mu_j=\exp(\alpha_j+{\rm i}\varphi_j).\]
	The spectrum of $D_{\Gamma(A)}$ is then given by
	\begin{equation*}
		\sigma\bigl(D_{\Gamma(A)}\bigr)=\{(\varphi_j-{\rm i}\alpha_j+2\pi k)/l \mid j\in\{1,\dots,m\}\text{, } k\in\mathbb{Z} \}.
	\end{equation*}
\end{Corollary}

\begin{Remark}
	Similar results exist for the Laplace operator on metric graphs, see for exam\-ple~\mbox{\cite[Theorem 2.1.8]{BerKuch}}, \cite[Theorem 19]{Kuchment2014} and the main theorems in \cite{VONBELOW}.
\end{Remark}

\begin{proof}
	By Theorem~\ref{thm:spec}, $\lambda$ is an eigenvalue of $D_{\Gamma(A)}$ if and only if ${\rm e}^{{\rm i}\lambda l}\in\mathbb{C}$ is an eigenvalue of~$A$. This proves the corollary.
\end{proof}

The following lemma allows us to delete vertices that have exactly one in-coming and one outgoing edge when dealing with the multi-dimensional characteristic polynomial. This will facilitate future arguments and calculations regarding the spectrum of the Dirac operator.

Consider the following setup: Let $G$ be a metric directed multigraph and assume that there exists a vertex $v\in V$ with in- and out-degree being one. Furthermore, assume that the ingoing edge at $v$ is different from the outgoing edge. Enumerate the edges of $G$ by $e_1,e_2,\dots,e_n$ such that $e_1$ is the ingoing and $e_2$ is the outgoing edge at $v$, defining a basis of $\mathbb{C}^E$. Let $\tilde{G}=\big(\tilde{V},\tilde{E}\big)$ be the graph obtained out of $G$ by replacing $e_1$ and $e_2$ with a single edge $e$ running from the tail of $e_1$ to the head of $e_2$. We define a length function $\tilde{l}\colon \tilde{E}\rightarrow (0,\infty)$ by
\begin{align*}
	&\tilde{l}_e=l_{e_1}+l_{e_2},\\
	&\tilde{l}_{e_i}= l_{e_i} \qquad \forall i>2.
\end{align*}
\begin{figure}[h]
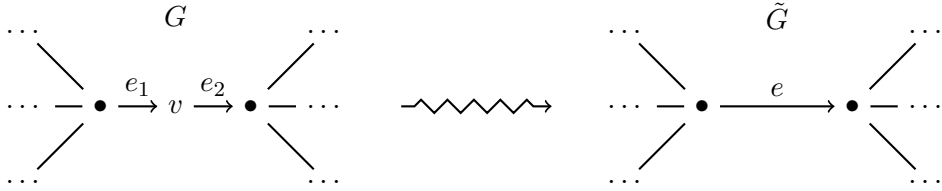

	\centering
	\tikz [thick] {
		\node (G) at (2,1.2) {$G$};
		\node (a1) at (0,1) {$\dots$};
		\node (b1) at (0,0) {$\dots$};
		\node (c1) at (0,-1) {$\dots$};
		\node (v1) at (1,0) {$\bullet$};
		\node (v2) at (2,0) {$v$};
		\node (v3) at (3,0) {$\bullet$};
		\node (d1) at (4,1){$\dots$};
		\node (e1) at (4,0){$\dots$};
		\node (f1) at (4,-1) {$\dots$};
		
		\graph { {(a1),(b1),(c1)} -- (v1) ->["$e_1$"] (v2) ->["$e_2$"] (v3) --{(d1),(e1),(f1)} };

		\draw[->, decorate, decoration={zigzag, pre length = 5pt, post length = 5pt}] (5,0) -- (7,0);
		
		\node (G2) at (10,1.2) {$\tilde{G}$};
		\node (a2) at (8,1) {$\dots$};
		\node (b2) at (8,0) {$\dots$};
		\node (c2) at (8,-1) {$\dots$};
		\node (w1) at (9,0) {$\bullet$};
		\node (w3) at (11,0) {$\bullet$};
		\node (d2) at (12,1){$\dots$};
		\node (e2) at (12,0){$\dots$};
		\node (f2) at (12,-1) {$\dots$};
		
		\graph { {(a2),(b2),(c2)} -- (w1) ->["$e$"] (w3) --{(d2),(e2),(f2)} };
	}
	\caption{Removal of a vertex.}
\end{figure}

\begin{Lemma}\label{lem:reduction}
Consider the above setup and let $A$ be a $G$-endomorphism. The transformation matrix of $A$ with respect to the basis $e_1,\dots,e_n$ of $\mathbb{C}^E$ is then of the form
	\[
	A=
		\begin{pmatrix}
			\begin{matrix}
				0 & \alpha_{12} \\
				\alpha & 0
			\end{matrix}
			&
			\begin{matrix}
				\ldots &\alpha_{1n}\\
				\ldots & 0
			\end{matrix}
			\\
			\begin{matrix}
				0 & \alpha_{32}\\
				\vdots & \vdots\\
				0 & \alpha_{n2}
			\end{matrix}
			&
			\hat{A}
		\end{pmatrix}.
	\]
	Let $\tilde{A}$ be the endomorphism on $\mathbb{C}^{\tilde{E}}$ that with respect to basis $e,e_3,\dots,e_n$ has the transformation matrix
	\[
	\tilde{A}=
	\begin{pmatrix}
		\alpha\alpha_{12}
		&
		\begin{matrix}
			\ldots &\alpha\alpha_{1n}
		\end{matrix}
		\\
		\begin{matrix}
			\alpha_{32}\\
			\vdots \\
			\alpha_{n2}
		\end{matrix}
		&
		\hat{A}
	\end{pmatrix}.
	\]
	Then $\tilde{A}$ is a $\tilde{G}$-endomorphism such that for the multi-dimensional characteristic polynomials
	\[
		P_A(x_1,x_2,\dots,x_n)=P_{\tilde{A}}(x_1 x_2,\dots,x_n)
	\]
	holds.
	In particular, we have
	\begin{equation*}
		P_{A,l} = P_{\tilde{A},\tilde{l}}.
	\end{equation*}
\end{Lemma}

\begin{proof}
	The lemma follows from a quick calculation using the Laplace expansion along the first column and the multilinearity of determinants:
	\begin{align*}
		P_A(x_1,x_2,\dots,x_n)&=\det
			\begin{pmatrix}
				\begin{matrix}
				\hphantom{-} x_1 & -\alpha_{12} \\
					-\alpha & \ x_2
				\end{matrix}
				&
				\begin{matrix}
					\ldots &-\alpha_{1n}\\
					\ldots & 0
				\end{matrix}
				\\
				\begin{matrix}
					\hphantom{-}0 & -\alpha_{32}\\
					\hphantom{-}\vdots & \vdots\\
					\hphantom{-}0 & -\alpha_{n2}
				\end{matrix}
				&
				\diag(\hat{x})-\hat{A}
			\end{pmatrix}\\
		&=x_1 \det
				\begin{pmatrix}
				\ x_2
				&
				\begin{matrix}
					0 & \ldots & 0
				\end{matrix}
				\\
				\begin{matrix}
					-\alpha_{32}\\
					\vdots\\
					-\alpha_{n2}
				\end{matrix}
				&
				\diag(\hat{x})-\hat{A}
			\end{pmatrix}
		+\alpha \det
			\begin{pmatrix}
				-\alpha_{12}
				&
				\begin{matrix}
					\ldots & -\alpha_{1n}
				\end{matrix}
				\\
				\begin{matrix}
					-\alpha_{32}\\
					\vdots\\
					-\alpha_{n2}
				\end{matrix}
				&
				\diag(\hat{x})-\hat{A}
			\end{pmatrix}\\
		&= \det
			\begin{pmatrix}
				x_1 x_2 -\alpha\alpha_{12}
				&
				\begin{matrix}
					\ldots & -\alpha\alpha_{1n}
				\end{matrix}
				\\
				\begin{matrix}
					-\alpha_{32}\\
					\vdots\\
					-\alpha_{n2}
				\end{matrix}
				&
				\diag(\hat{x})-\hat{A}
			\end{pmatrix}=P_{\tilde{A}}(x_1 x_2, \dots,x_n).
	\end{align*}
	It is easily checked that $\tilde{A}$ indeed defines a $\tilde{G}$-automorphism. The equality of the characteristic functions follows from
	\[
		P_{A,l}(\lambda)=P_A\big({\rm e}^{{\rm i}\lambda l_1},{\rm e}^{{\rm i}\lambda l_2},\dots,{\rm e}^{{\rm i}\lambda l_n}\big)=P_{\tilde{A}}\big({\rm e}^{{\rm i}\lambda (l_1+l_2)},\dots,{\rm e}^{{\rm i}\lambda l_n}\big)=P_{\tilde{A},\tilde{l}}(\lambda).
\tag*{\qed}
\]
\renewcommand{\qed}{}
\end{proof}

\begin{Example}\label{exa:cycle}
	Let $C_n$ be the directed cycle with $n$ edges and let $A$ be a unitary $C_n$-en\-do\-mor\-phism. Enumerating the edges $e_1,\dots,e_n$ in cyclic order, it is easily seen that there exists $\varphi_1,\dots,\varphi_n\in[0,2\pi)$ such that the transformation matrix with respect to this basis is of the form
	\[
	A=
	\begin{pmatrix}
		0& & &{\rm e}^{{\rm i}\varphi_n}\\
		{\rm e}^{{\rm i}\varphi_1}&0 & & \\
		& &\ddots &\\
		& & {\rm e}^{{\rm i}\varphi_{n-1}} &0
	\end{pmatrix}.
	\]
	By removing all vertices but one, according to Lemma \ref{lem:reduction}, the multi-dimensional characteristic polynomial is of the form
	\begin{align*}
		P_A(x_1,\dots,x_n) =x_1\cdots x_n-\exp({\rm i}\varphi_1)\cdots \exp({\rm i}\varphi_n) =\prod_{i=1}^n x_i - \exp\Biggl({\rm i}\sum_{i=1}^n\varphi_i\Biggr).
	\end{align*}
	Hence, the characteristic function is given by
	\begin{align*}
		P_{A,l}(\lambda)
		= \exp\Biggl({\rm i}\lambda\sum_{j=1}^n l_j\Biggr)-\exp\Biggl({\rm i}\sum_{j=1}^n\varphi_j\Biggr).
	\end{align*}
	This implies that the spectrum of the canonical Dirac operator with respect to the boundary condition $\Gamma(A)$ is given by
	\begin{equation*}
		\lambda_k=\frac{1}{L}\Biggl(\sum_{j=1}^n\varphi_j +2\pi k\Biggr) \qquad \text{for}\ k\in\mathbb{Z},
	\end{equation*}
	where $L=\sum_{j=1}^n l_j$. By deleting the last column of $ \diag\big({\rm e}^{{\rm i}\lambda_k l}\big)-A $, we obtain a matrix of the form
	\begin{equation*}
		\tilde{A}_k=
		\begin{pmatrix}
			{\rm e}^{{\rm i}\lambda_k l_1} & & &\\
			-{\rm e}^{{\rm i}\varphi_1}& {\rm e}^{{\rm i}\lambda_k l_2} & &\\
			&-{\rm e}^{{\rm i}\varphi_2} &\ddots& \\
			& & & {\rm e}^{{\rm i}\lambda_k l_{n-1}}\\
			& & & -{\rm e}^{{\rm i}\varphi_{n-1}}
		\end{pmatrix}.
	\end{equation*}
	Clearly, $ \tilde{A}_k $ has rank $n-1$. Thus, the kernel of $ \diag\bigl({\rm e}^{{\rm i}\lambda_k l}\bigr)-A $ can only be one-dimensional, which implies that the multiplicity of $\lambda_k$ is one.
\end{Example}

In the following, we will be dealing with multiple graphs at once, so it is necessary to establish a notation for this situation. Given a graph $G$, we denote the set of edges of $G$ by $E(G)$ and the set of vertices of $G$ by $V(G)$.

A matrix $A$ is reducible if there exists a permutation matrix $P$ such that $P^{-1}AP$ is of the form
\[
	\begin{pmatrix}
		X & 0\\
		Y & Z\\
	\end{pmatrix},
\]
where $X$ and $Z$ are square matrices. Otherwise, $A$ is said to be irreducible \cite[p.~18]{Cvet}.
We transfer this definition to linear maps on $G$.

\begin{Definition}
	A $G$-endomorphism $A$ is said to be {\it reducible} with respect to a proper subset~$\Delta$ of $E(G)$ if $A$ maps $\mathbb{C}^{\Delta}:=\bigl\{\sum_{e\in \Delta}c_e \cdot e \mid c_e\in \mathbb{C}\bigr\}$ into itself. If no such $\Delta$ exists, $A$ is said to be {\it irreducible}.
\end{Definition}

In the case that $G$ is disconnected and $G^\prime$ is a connected component of $G$, then all $G$-endomorphisms $A$ are reducible with respect to $E(G^\prime)$.

Given $\Delta\subseteq E(G)$, we denote its complement by $\Delta^c$ and by $G[\Delta]$ the subgraph of $G$ induced by $\Delta$. By the restriction of a $G$-endomorphism $A$ to $\mathbb{C}^{\Delta}$ we mean the restriction of
\[
	\pi_{\Delta}\circ A
\]
to $\mathbb{C}^{\Delta}$, where $\pi_{\Delta}\colon \mathbb{C}^{E}\rightarrow \mathbb{C}^{\Delta}$ is the canonical projection map.

\begin{Theorem}\label{prop:decomp}
	Let $A$ be a $G$-endomorphism. If $A$ is reducible with respect to $\Delta\subsetneq E(G)$, then the spectrum of $D_{\Gamma(A)}$ is given by
	\begin{equation*}
		\sigma\bigl(D_{\Gamma(A)}\bigr)=\sigma\big(D_{\Gamma(A_1)}^1\big)\cup\sigma\big(D_{\Gamma(A_2)}^2\big),
	\end{equation*}
	where $A_1$ is the restriction of $A$ to $ \mathbb{C}^{\Delta} $ and $A_2$ is the restriction of $A$ to $ \mathbb{C}^{\Delta^c} $. The operators $D^1$ and $D^2$ are the Dirac operators on $G[\Delta]$ and $G[\Delta^c]$, respectively.
\end{Theorem}

\begin{proof}
	Enumerate the edges $e_1,\dots,e_n$ of $G$ such that $\{e_1,\dots,e_k\}=\Delta$. Then, the transformation matrix of $A$ with respect to this basis is of the form
	\[
		\begin{pmatrix}
			A_1 & C\\
			0 & A_2\\
		\end{pmatrix},
	\]
	where $A_1$ and $A_2$ are square matrices. The matrix $A_1$ represents the restriction of $A$ to $\mathbb{C}^{\Delta}$ and~$A_2$ represents~$A$ restricted to $ \mathbb{C}^{\Delta^c} $. These matrices give boundary conditions on $G[\Delta]$ and $G[\Delta^c]$, respectively. Denote the induced length function on $G[\Delta]$ and $G[\Delta^c]$ by $l_1$ and $l_2$, respectively. We then have
	\begin{align*}
		P_{A,l}(\lambda)=
		\det
		\begin{pmatrix}
			\diag(\exp({\rm i}\lambda l_1))-A_1 & -C\\
			0 & \diag(\exp({\rm i}\lambda l_2))-A_2\\
		\end{pmatrix}
		=P_{A_1,l_1}(\lambda)\cdot P_{A_2,l_2}(\lambda).
	\end{align*}
	Therefore, $P_{A,l}(\lambda)$ is zero if either $P_{A_1,l_1}(\lambda)$ or
	$P_{A_2,l_2}(\lambda)$ is zero and, hence,
	\[
		\sigma(D_A)=\sigma\big(D_{\Gamma(A_1)}^1\big)\cup\sigma\big(D_{\Gamma (A_2)}^2\big).
\tag*{\qed}
\]
\renewcommand{\qed}{}
\end{proof}

\begin{Example}
	Let $G$ be a directed multigraph such that there exists a partition $\{V_1, V_2\}$ of the vertex set $V$ of $G$ so that there is no edge running from $V_1$ to $V_2$. For $i,j \in \{1,2\}$, define
	\begin{equation*}
		E_{ij}= \bigl\{e\in E\mid \partial_- e\in V_i, \partial_+ e\in V_j\bigr\}.
	\end{equation*}
	By the assumptions on $G$ we have $E_{12}=\varnothing$ and a partition
	\begin{equation*}
		E = E_{11}\cup E_{21}\cup E_{22}.
	\end{equation*}
	Consider the subgraphs
	\begin{alignat*}{3}
		&G_1= G[E_{11}], && G_2= G[E_{21}\cup E_{22}],&\\
		&G_{21}= G[E_{21}],\qquad && G_{22}= G[E_{22}].&
	\end{alignat*}
	Then every $G$-endomorphism $A$ is reducible with respect to $E_{11}$. If $A_{2}$ denotes the restriction of $A$ to $\mathbb{C}^{E_{21}\cup E_{22}}$, we also have that $A_2$ is reducible with respect to $E_{21}$ with the restriction $A_{21}$ to $\mathbb{C}^{E_{21}}$ being zero. According to Theorem~\ref{prop:decomp} and \ref{thm:spec}, we obtain for the spectrum
	\begin{align*}
		\sigma\bigl(D_{\Gamma(A)}\bigr)&= \sigma\big(D^1_{\Gamma(A_1)}\big)\cup \underbrace{\sigma\big(D^{21}_{\Gamma(A_{21})}\big)}_{=\varnothing} \cup \sigma\big(D^{22}_{\Gamma(A_{22})}\big)\\
		&= \sigma\big(D^1_{\Gamma(A_1)}\big) \cup \sigma\big(D^{22}_{\Gamma(A_{22})}\big),
	\end{align*}
	where $D^{1}$, $D^{21}$, $D^{22}$ are the scalar Dirac operators of $G_{1}$, $G_{21}$, $G_{22}$ and $A_1$, $A_{21}$, $A_{22}$ are the restrictions to $\mathbb{C}^{E_1}$, $\mathbb{C}^{E_{21}}$, $\mathbb{C}^{E_{22}}$ respectively.
	
	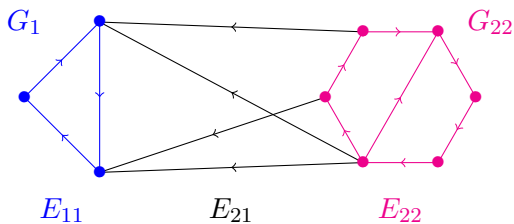
\begin{figure}[t]
		\centering
		\begin{tikzpicture}
			[decoration={
				markings,
				mark=at position 0.5 with {\arrow{>}}}]

			\draw[postaction={decorate}] (4.5,-0.87)--(1,-1);
			\draw[postaction={decorate}] (4,0)--(1,-1);
			\draw[postaction={decorate}] (4.5,-0.87)--(1,1);
			\draw[postaction={decorate}] (4.5,0.87)--(1,1);
			
			\node[blue] (G_1) at (0,1) {$G_1$};
			\filldraw[blue] (0,0) circle (1.5pt) node {$\bullet$};
			\filldraw[blue] (1,1) circle (1.5pt) node {$\bullet$};
			\filldraw[blue] (1,-1) circle (1.5pt) node {$\bullet$};
			\draw[blue, postaction={decorate}] (0,0)--(1,1);
			\draw[blue, postaction={decorate}] (1,1)--(1,-1);
			\draw[blue, postaction={decorate}] (1,-1)--(0,0);
			
			\node[magenta] (G_2) at (6.2,1) {$G_{22}$};
			\filldraw[magenta] (4,0) circle (1.5pt) node {$\bullet$};
			\filldraw[magenta] (4.5,0.87) circle (1.5pt) node {$\bullet$};
			\filldraw[magenta] (5.5,0.87) circle (1.5pt) node {$\bullet$};
			\filldraw[magenta] (6,0) circle (1.5pt) node {$\bullet$};
			\filldraw[magenta] (5.5,-0.87) circle (1.5pt) node {$\bullet$};
			\filldraw[magenta] (4.5,-0.87) circle (1.5pt) node {$\bullet$};
			\draw[magenta, postaction={decorate}] (4,0)--(4.5,0.87);
			\draw[magenta, postaction={decorate}] (4.5,0.87)--(5.5,0.87);
			\draw[magenta, postaction={decorate}] (5.5,0.87)--(6,0);
			\draw[magenta, postaction={decorate}] (6,0)--(5.5,-0.87);
			\draw[magenta, postaction={decorate}] (5.5,-0.87)--(4.5,-0.87);
			\draw[magenta, postaction={decorate}] (4.5,-0.87)--(5.5,0.87);
			\draw[magenta, postaction={decorate}] (4.5,-0.87)--(4,0);

			\node[blue] (E_1) at (0.5,-1.5) {$E_{11}$};
			\node (E_2) at (2.75,-1.5) {$E_{21}$};
			\node[magenta] (E_3) at (5,-1.5) {$E_{22}$};
		\end{tikzpicture}
		\caption{An example of a graph that has two subgraphs $G_1$, $G_{22}$ with no edge going from $G_1$ to $G_{22}$.}
	\end{figure}
	
\end{Example}

\section{Boundary conditions}

In this section, we will look at special boundary conditions for the scalar Dirac operator.

\subsection{Decompositions of Eulerian directed multigraphs}

Recall that a directed walk is a sequence of edges
\[
	W:=e_1,e_2,\dots, e_k,
\]
where $\partial_+ e_{i-1}= \partial_- e_i$. Denote by $E(W)$ the set of edges appearing in $W$. A directed walk is said to be closed if $ \partial_+ e_k= \partial_- e_0$. A trail is a walk in which all edges are distinct. A closed directed trail that does not pass a vertex twice, except for $\partial_- e_0$, is called a cycle. A cycle/trail decomposition of $G$ is a family
\[
	S=\{C_1,\dots,C_m\}
\]
of cycles/trails such that every edge of $G$ appears in a unique cycle/trail of $S$.

\begin{Definition}
	Let $G$ be a directed multigraph. We call a bijective map $P\colon E(G)\rightarrow E(G)$ a~{\it permutation of the edges} of $G$, the set of permutations of its edges are denoted by $\mathfrak{S}(E(G))$. If $P$ satisfies that $P(e)=e^\prime$ implies $\partial_+e= \partial_- e^\prime$, then~$P$ is called a {\it $G$-permutation}. We denote the set of all $G$-permutations by $\mathfrak{S}(G)$.
\end{Definition}

\begin{Remark}
	We have $\mathfrak{S}(G)\subseteq \mathfrak{S}(E(G))$. By enumerating $E(G)$, we assign to every permutation $P\in \mathfrak{S}(E(G))$ a corresponding element in the symmetric group $\mathfrak{S}_n$, where $n$ is the number of edges of $G$.
\end{Remark}

\begin{Lemma}
	A permutation $P\in \mathfrak{S}(G)$ of $G$ induces for each vertex $v$ a bijection $P_v$ from the ingoing edges of $v$ to the outgoing edges of $v$. In particular, $\mathfrak{S}(G)$ is empty if $G$ is not a Eulerian directed multigraph.
\end{Lemma}

\begin{proof}
	Fix $v\in V(G)$. Let $e$ be an ingoing vertex at $v$. By the definition of $G$-permutations, $P(e)$ is then an outgoing vertex of $v$. Therefore, $P_v$ is well-defined. Injectivity of $P$ implies injectivity of $P_v$ and, since $P$ is surjective, for every outgoing edge~$e^\prime$ at $v$, there exists an edge~$e$ of~$G$ mapped by $P$ onto $e^\prime$. This implies that the head of $e$ is the tail of $e^\prime$. Therefore, $e$ is an ingoing edge of $v$ and $P_v$ maps $e$ onto $e^\prime$, which proves that $P_v$ is bijective.
\end{proof}

It is clear that every $G$-permutation $P\in\mathfrak{S}(G)$ induces a unitary $G$-endomorphism, and hence, a local self adjoint boundary condition. We will now investigate the spectrum of the Dirac operator subject to such a boundary condition.

To this end, recall the following facts concerning permutations: for every permutation $\pi\in \mathfrak{S}_n$ and $x\in\{1,\dots,n\}$ there exists a $k\in \mathbb{N}$ such that $\pi^k(x)=x$. If $k$ is the smallest such integer, then $x,\pi(x),\pi^2(x),\dots,\pi^{k-1}(x)$ is said to be a $k$-cycle in $\pi$. Every permutation can be decomposed into a disjoint union of its cycles.

\begin{Theorem}\label{thm:perm}
	Let $G$ be a Eulerian directed multigraph. Every decomposition of the directed multigraph into directed closed trails corresponds to a $G$-permutation $P$ and vice versa. Each trail of length $k$ then corresponds to a $k$-cycle of $P$, i.e., the trail containing the edge $e$ is given by \looseness=-1
	\begin{equation*}
		e,P(e),P^2(e),\dots,P^{k-1}(e).
	\end{equation*}
\end{Theorem}

\begin{proof}
	Given a decomposition of the directed multigraph into closed directed trails, define a~map $P\colon E(G)\rightarrow E(G)$ by $P(e)=e^\prime$, where $e^\prime$ is the edge that comes after $e$ in the closed trail containing $e$. In order to show that $P$ is a permutation, it remains to be shown that it is injective. Hence, assume that
	\begin{equation*}
		P(e)=e^\prime=P(\tilde{e})
	\end{equation*}
	and, therefore, there is a trail $C$ that contains $e$ and a trail $\tilde{C}$ that contains $\tilde{e}$ such that $e^\prime$ is the subsequent edge in both trails. Hence, $C$ and $\tilde{C}$ both contain $e^\prime$ and, since the decomposition demands the trails be disjoint, we get $C=\tilde{C}$. A trail cannot contain an edge twice and, therefore, $e$ is equal to $\tilde{e}$. Consequently, injectivity of $P$ follows.
	
	Now, let $P$ be a $G$-permutation. For an edge $e$, define a directed walk
	\begin{equation}\label{eq:trail}
		T(e)=e,P(e),P^2(e),\dots,P^{k-1}(e),
	\end{equation}
	where $k$ is the minimum number such that $P^{k}(e)=e$. The fact that $P$ is a permutation of the set of edges of $G$ ensures that such a $k$ exists. It is, indeed, a directed walk since, by definition of $P$, the head of $P^i(e)$ is the tail of $P^{i+1}(e)$. It is closed, since the head of $P^{k-1}(e)$ is the tail of $P^{k}(e)=e$. It remains to be shown that it is a trail, i.e., that there exists no $l< j$ in $\{0,\dots,k-1\}$ such that $P^l(e)=P^j(e)$. This would imply $P^{j-l}(e)=e$, which would contradict the definition of $k$, since $j-l<k$.
	
	Hence, equation \eqref{eq:trail} defines, for every edge, a directed closed trail that contains this edge. To show that these trails are a decomposition of the directed multigraph, it remains to be proven that two distinct trails do not share any edges. So, let
	\begin{align*}
		T_1=e,P(e),P^2(e),\dots,P^{k_1}(e),\qquad
		T_2=e^\prime,P(e^\prime),P^2(e^\prime),\dots,P^{k_2}(e^\prime)
	\end{align*}
	be two trails that share an edge, i.e., $P^{m_1}(e)=P^{m_2}(e^\prime)$ for some $1\leq m_i\leq k_i$. Without loss of generality, let $m_1\geq m_2$. We then arrive at $P^{m_1-m_2}(e)=e^\prime$. Applying $P^{k_1}$, which commutes with $P^{m_1-m_2}$, to both sides of this equation, we see that
	\[
		P^{m_1-m_2}(e)=P^{k_1}(e^\prime).
	\]
	Therefore, $P^{k_1}(e^\prime)$ is equal to $e^\prime$, which implies that $k_1\geq k_2$. The inequality $k_2\geq k_1$ follows analogously, which shows that $k_1=k_2$.
	Therefore, $P^{m_1-m_2}$ gives a bijection between the edges of $T_1$ and $T_2$. Therefore, they are the same trail.
	
	It is clear that these constructions are inverse to each other and, hence, this correspondence is one to one.
\end{proof}

\begin{Proposition}
	Let $P$ be a $G$-permutation. Then, the trace of $P$ is equal to the number of loops in the decomposition of $G$ that corresponds to $P$.
\end{Proposition}
\begin{proof}
	The trace of $P$ is equal to
	\begin{equation*}
		\tr P = \sum_{e\in E(G)} \langle e, P(e)\rangle
	\end{equation*}
	and further
	\[\langle e, P(e)\rangle = \begin{cases}
		1 & \text{if}\ P(e)=e, \\ 0 & \text{else}.
	\end{cases}\]
	Hence, the trace of $P$ is equal to the number of edges such that $P(e)=e$. For this to hold, $e$ must be a loop, and furthermore by the proof of Theorem~\ref{thm:perm}, it is the closed trail containing~$e$.
\end{proof}

Given a closed trail $T=e_1,\dots,e_k$, let its length be
\[
	L_T=\sum_{i=1}^k l_{e_i}.
\]

\begin{Lemma}\label{lem:cycle}
	Let $P$ be a $G$-permutation. Then, for each closed trail $T$ in the decomposition defined by $P$, $P$ is reducible with respect to $E(T)$.
\end{Lemma}
\begin{proof}
	Let $T$ be a closed trail in the decomposition defined by $P$ and let $e$ be an edge that appears in $T$. Then, by definition of $T$, $P(e)$ is again an edge that appears in $T$, which shows that $P$ is reducible with respect to $E(T)$.
\end{proof}

We say that $\alpha\in \mathbb{R}$ is divisible by $\beta\in \mathbb{R}$ if there exists a $k\in\mathbb{Z}$ such that $ \alpha=k\cdot\beta $.

\begin{Theorem}\label{thm:perm2}
	Let $G$ be a Eulerian directed multigraph and let $P$ be a linear $G$-permutation. If $T_1,\dots,T_m$ are the closed trails of the decomposition defined by $P$, then the spectrum of $D_{\Gamma(P)}$ is given by
	\begin{equation*}
		\lambda_{k,i} = \frac{2\pi k}{L_i}\qquad \text{with} \ i\in\{1,\dots,m\}, \ k\in \mathbb{Z},
	\end{equation*}
	where $L_j$ is the length of $T_j$. The multiplicity $m(\lambda_{k,i})$ is, for $k\neq 0$, equal to the number of cycles whose length is divisible by $L_i/k$. Furthermore, $m(0)$ is equal to the total number of closed trails in the decomposition defined by $P$.
\end{Theorem}
\begin{Remark}
	If we choose all edges to have unit length, then $L_j$ is the actual number of edges appearing in the closed trail $T_j$.
\end{Remark}

\begin{proof}
	By Lemma \ref{lem:cycle}, $P$ is reducible with respect to each $T_i$ and, hence, by Theorem~\ref{prop:decomp}, we have that
	\begin{equation*}
		\sigma\big(D_{\Gamma(P)}\big)=\bigcup_{i=1}^m\sigma\big(D_{\Gamma(P_i)}\big),
	\end{equation*}
	where $P_i$ is the restriction of $P$ to $E(T_i)$. Hence we can assume that we have only one trail $T$.
	
	Enumerate the edges of $T$ in such a way that $P(e_i)= e_{i +1 \mod \lvert E(T)\rvert }$.
	The transformation matrix of $P$ with respect to the basis $e_{0},\dots,e_{\lvert E(T)\rvert-1}$ is equal to
	\begin{equation*}
		\begin{pmatrix}
			0 & & & &1\\
			1& 0 & & &\\
			& 1 & 0 & &\\
			& &\ddots &\ddots &\\
			& & & 1 & 0
		\end{pmatrix}.
	\end{equation*}
	As in Example \ref{exa:cycle}, the spectrum of $P$ can be easily calculated to be
	\begin{equation*}
		\sigma(D_{\Gamma(P)})=\left\{\frac{2\pi k}{L}\mid k\in \mathbb{Z}\right\},
	\end{equation*}
	with each eigenvalue having multiplicity one.
\end{proof}

\begin{Corollary}
	Let $P$ be a $G$-permutation. Let $\lambda_0$ be the smallest positive eigenvalue of $D_{\Gamma(P)}$. Then,
$L=\frac{2\pi}{\lambda_0}$
	is the length of the longest cycle in the decomposition associated with $P$. The multiplicity of~$\lambda_0$ is equal to the number of cycles of length $L$.
\end{Corollary}
\begin{proof}
	Let $\lambda_0$ be the smallest positive eigenvalue of $D_{\Gamma(P)}$. By Theorem~\ref{thm:perm2}, it is of the form
$\lambda_0= \frac{2\pi k}{L}$,
	where $k$ is an integer and $L$ is the length of a closed trail appearing in the decomposition defined by $P$. In order for $\lambda_0$ to be the smallest positive eigenvalue, $k$ has to be one and $L$ must be the maximum of the lengths of the closed trails.
\end{proof}

\subsection{The directed adjacency matrix}

As we have seen, the spectrum of the Dirac operator is closely tied to the boundary condition it is subjected to. Boundary conditions involve choices which are a priori not canonical. In order to obtain a theory that ties the spectrum of the Dirac operator to properties of the underlying graph, we need to come up with an intrinsic boundary condition. By Proposition \ref{prop: s.a. implies eulerian}, this boundary condition will in general not be self adjoint. However, canonical boundary conditions on graphs are given by directed adjacency matrices:

\begin{Definition}
	The {\it directed adjacency matrix $\mathcal{A}(G)$} of a directed multigraph $G$ is given by
	\begin{align*}
		&\mathcal{A}(G)= \sum_{e, e^\prime \in E} I_{e, e^\prime} e\otimes (e^\prime)^\ast \in \mathbb{C}^E\otimes\big(\mathbb{C}^E\big)^\ast = \mathrm{End}\bigl(\mathbb{C}^E\bigr), \\
		&I_{e, e^\prime}= \begin{cases}1 & \text{if} \ \partial_+ e^\prime = \partial_- e,\\ 0 & \text{otherwise}.\end{cases}
	\end{align*}
\end{Definition}

\begin{Lemma}\label{lem:inciauto}
	The directed adjacency matrix $\mathcal{A}(G)$ is a $G$-endomorphism.
\end{Lemma}
\begin{proof}
	Let $v\in V$ be a vertex of $G$ and $e^\prime$ be an edge in $\mathbb{C}^{E_v^+}$. Therefore, $\partial_+ e^\prime = v$, and
	\[\mathcal{A}(G)e^\prime=\sum_{e\in E(G)}I_{e,e^\prime}e=\sum_{\substack{e \\ \partial_- e= v}}e \in \mathbb{C}^{E_v^-},\]
	which proves the lemma, since $\mathbb{C}^{E_v^+}$ is spanned by $\{ e^\prime \mid \partial_+ e^\prime = v\}$.
\end{proof}

Lemma \ref{lem:inciauto} shows that the adjacency matrix defines a boundary condition for the Dirac operator. We denote the Dirac operator subject to the boundary condition $\Gamma(\mathcal{A}(G))$ by $D_G$. As we will see, the spectrum of $D_G$ highly depends on the directed cycles present in $G$. We use the notation $G^\prime\subseteq G$ to say that $G^\prime$ is a subgraph of $G$ with no isolated vertices.

\begin{Definition}
	A {\it disjoint collection of cycles} is a directed multigraph $C$, whose connected components consist of directed cycles. The {\it connectedness $\alpha(C)$}
of $C$ is the number of connected components of $C$, the {\it size $\eta(C)$}
of $C$ is the number of edges of $C$. For a directed multigraph $G$, set
	\[
	\mathcal{C}(G):=\{C\subseteq G\mid C \text{ is a disjoint collection of cycles} \}.
	\]
\end{Definition}

\begin{Remark}The empty subgraph $\varnothing$ with no vertices and edges is always an element of $\mathcal{C}(G)$. Its connectedness and size equal zero.
\end{Remark}

The following lemma gives a characterization of disjoint collections of cycles:

\begin{Lemma}\label{lem:idod}
	If $G$ is a directed multigraph such that for every vertex $v\in V(G)$ the in-degree~${\rm id}(v)$ and the out-degree ${\rm od}(v)$ are both equal to one, then $G$ is a disjoint collection of cycles.
\end{Lemma}

\begin{proof}It is easily checked that, in this case, $\mathcal{A}(G)$ is a $G$-permutation. By Theorem~\ref{thm:perm}, the directed adjacency matrix defines a decomposition $S$ of $G$ into directed trails. However, if a trail in $S$ passes a vertex twice or two distinct trails intersect at a vertex, then the in- and out-degrees of the vertex would be greater than one, which would contradict the assumption. Consequently, $G$ is a disjoint collection of cycles.
\end{proof}

\begin{Proposition}\label{prop:adjacency}
	The directed adjacency matrix $\mathcal{A}(G)$ of $G$ is non-singular if and only if $G$ is a disjoint collection of directed cycles. In this case, the determinant of $\mathcal{A}(G)$ is given by the formula
	\[
		\det \mathcal{A}(G) = (-1)^{\alpha(G)-\eta(G)}.
	\]
\end{Proposition}

\begin{proof}
	Let $v$ be a vertex of $G$. If the in-degree/out-degree of $v$ is equal to zero, then the out-degree/in-degree must be unequal zero as $G$ has no isolated vertices. Then, for every edge~$e$ starting/ending in $v$, the row/column of $\mathcal{A}(G)$ that corresponds to $e$ is equal to zero, which implies that $\mathcal{A}(G)$ is singular. Therefore, in order for $\mathcal{A}(G)$ to be non-singular, the in-degree and the out-degree of every vertex of $G$ must be greater or equal to one.
	
	Now, assume that the in-degree/out-degree of a vertex $v$ is greater than or equal to two and let $e_1$ and $e_2$ be two edges ending/starting in $v$. Then, the columns/rows of $\mathcal{A}(G)$ that correspond to $e_1$ and $e_2$ are equal, which implies that $\mathcal{A}(G)$ is singular. Consequently, $\mathcal{A}(G)$ can only be non-singular if the in-degree and the out-degree of all vertices are equal to one.
	
	Lemma \ref{lem:idod} implies that $G$ is a disjoint collection of cycles and $\mathcal{A}(G)$ is a $G$-permutation, where each cycle of $G$ corresponds to a cycle of $\mathcal{A}(G)$ of the same size. Writing $G$ as the disjoint union of cycles
	\[
		G=C_{k_1}\sqcup\dots \sqcup C_{k_s},
	\]
	where $\alpha(G)=s$ and $\eta(G)=k_1+\dots +k_s$, it follows that $\mathcal{A}(G)$ is the product of cycles of size~$k_1,\dots,k_s$. Since the signature of a cycle of size $k$ is given by $(-1)^{k-1}$, the signature $\sgn(\mathcal{A}(G))$ of $\mathcal{A}(G)$ is given by $-1$ to the power of
	\[
		\sum_{i=1}^s (k_i -1)=\eta(G)-\alpha(G).
	\]
	The determinant of a permutation matrix is equal to the signature of the underlying permutation, which implies that
	\[
		\det \mathcal{A}(G) = \sgn(\mathcal{A}(G)) = (-1)^{\eta(G)-\alpha(G)}=(-1)^{\alpha(G)-\eta(G)},
	\]
	and proves the proposition.
\end{proof}

For a subgraph $H$ of $G$, define a monomial in the polynomial ring over the variables $\{x_e\mid e\in E(G)\}$ by
\[
	x^H:=\prod_{e\not\in E(H)}x_e.
\]
It seems counter-intuitive to multiply over the edges that are not contained in $E(H)$, however, it is the most convenient notation for the next theorem.

\begin{Theorem}
	The multi-dimensional characteristic polynomial $P_{\mathcal{A}(G)}$ of the directed adjacency matrix $\mathcal{A}(G)$ is given by
	\begin{equation}\label{eq:charpol}
		P_{\mathcal{A}(G)}(x)=\sum_{C\in \mathcal{C}(G)}(-1)^{\alpha(C)}x^C.
	\end{equation}
\end{Theorem}

\begin{proof}
	First, assume that $G$ has no loops. The multi-dimensional characteristic polynomial is given by
	\begin{align*}
		P_{\mathcal{A}(G)}(x) = \det(\diag(x)-\mathcal{A}(G))
		 = \sum_{\sigma\in \mathfrak{S}(E(G))}\sgn(\sigma)\prod_{e\in E(G)}(\diag(x)-\mathcal{A}(G))_{\sigma(e),e}.
	\end{align*}
	For a subgraph $H\subseteq G$, let $\mathfrak{S}_H$ be the subset of $\mathfrak{S}(E(G))$, which consists of all permutation whose fixed point set is equal to $E(G)\setminus E(H)$. We can regard $\mathfrak{S}_H$ as the subset of $\mathfrak{S}(E(H))$ of permutations without fixed points. We can rewrite $P_{\mathcal{A}(G)}$ as
	\[
 P_{\mathcal{A}(G)}(x)=\sum_{H\subseteq G}\prod_{e\notin E(H)}(x_e-I_{e,e})\sum_{\sigma\in \mathfrak{S}_H}\sgn(\sigma)\prod_{e\in E(H)}(-I_{\sigma(e),e}),
 \]
	where we sum over all induced subgraphs of $G$. Since $G$ has no loops, $I_{e,e}$ is always zero, which simplifies the above formula to
	\[
 P_{\mathcal{A}(G)}(x)=\sum_{H\subseteq G}(-1)^{|E(H)|}x^H\sum_{\sigma\in \mathfrak{S}_H}\sgn(\sigma)\prod_{e\in E(H)}\mathcal{A}(G)_{\sigma(e),e}.
 \]
	Moreover, whenever 
$\sigma\in \mathfrak{S}(E(H))$ has a fixed point, the product $\prod_{e\in E(H)}\mathcal{A}(G)_{\sigma(e),e}$ is equal to zero. Therefore, we can sum over all $\sigma\in \mathfrak{S}(E(H))$ instead of $\mathfrak{S}_H$. Further, the restriction of~$\mathcal{A}(G)$ to the edges of $H$ gives us the directed adjacency matrix of $H$, i.e.,
	\[\mathcal{A}(H)=\sum_{e,e^\prime\in E(H)}I_{e,e^\prime}e\otimes e^\prime.\]
	Altogether, this implies that
	\[P_{\mathcal{A}(G)}=\sum_{H\subseteq G}(-1)^{|E(H)|}x^H\det \mathcal{A}(H).\]
	Using Proposition \ref{prop:adjacency}, the multi-dimensional characteristic polynomial takes the form{\samepage
	\[P_{\mathcal{A}(G)}=\sum_{C\in\mathcal{C}(G)}(-1)^{\eta(C)}x^C(-1)^{\alpha(C)-\eta(C)}=\sum_{C\in\mathcal{C}(G)}(-1)^{\alpha(C)}x^C,\]
	which proves the theorem for the case in which $G$ has no loops.}
	
	Now, assume that $G$ has exactly one loop $e$. Split the loop into two edges $e_1$, $e_2$ by adding a~vertex $v$ to the loop to obtain a loopless graph $H$, see Figure \ref{figure:vertexadded}. We want to use Lemma~\ref{lem:reduction}: note that the graph $\tilde{H}$ obtained by removing $v$ is equal to $G$ and, moreover, the reduction of~$\mathcal{A}(H)$ to $G$ is equal to the directed adjacency matrix of $G$. We therefore have that
	\[P_{\mathcal{A}(G)}(x_e,x_3,\dots,x_n)=P_{I(H)}(x_e,1,x_3,\dots,x_n)=\sum_{\hat{C}\in\mathcal{C}(H)}(-1)^{\alpha(\hat{C})}(x_e,1,x_3,\dots,x_n)^{\hat{C}}.\]
	By replacing the cycle $e_1$, $e_2$ with $e$, it is clear that every collection of disjoint cycles $\hat{C}\in\mathcal{C}(H)$ corresponds to a unique collection $C\in\mathcal{C}(G)$ of same connectedness such that $e_1$, $e_2$ is a cycle in~$\hat{C}$ if and only if $e$ is a cycle in $C$.
	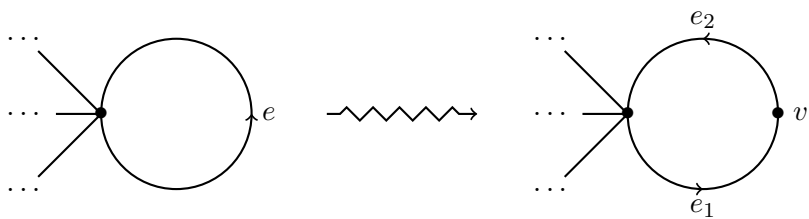
\begin{figure}[H]
		\centering
		\begin{tikzpicture}
			[thick, decoration={
				markings,
				mark=at position 0.5 with {\arrow{>}}}]
			\node (a1) at (0,1) {$\dots$};
			\node (b1) at (0,0) {$\dots$};
			\node (c1) at (0,-1) {$\dots$};
			\node (v1) at (1,0) {$\bullet$};

			\draw (a1) -- (1,0);
			\draw (b1) -- (1,0);
			\draw (c1) -- (1,0);
			\draw[postaction={decorate}] (v1)arc(180:540:1) node[pos=0.5, anchor=west] {$e$};

			\draw[->, decorate, decoration={zigzag, pre length = 5pt, post length = 5pt}] (4,0) -- (6,0);
			
			\node (a2) at (7,1) {$\dots$};
			\node (b2) at (7,0) {$\dots$};
			\node (c2) at (7,-1) {$\dots$};
			\node (v2) at (8,0) {$\bullet$};
			\node (v) at (10,0) {$\bullet$};
			\node (w) at (10.3,0) {$v$};
			
			\draw (a2) -- (8,0);
			\draw (b2) -- (8,0);
			\draw (c2) -- (8,0);
			\draw[postaction={decorate}] (v2)arc(180:360:1) node[pos=0.5, anchor=north] {$e_1$};
			\draw[postaction={decorate}] (v)arc(0:180:1) node[pos=0.5, anchor=south] {$e_2$};
		\end{tikzpicture}
		\caption{Adding a vertex $v$ to a loop $e$.}
		\label{figure:vertexadded}
	\end{figure}
	
	Then, the monomial $(x_e,x_3,\dots,x_n)^C$ is equal to $(x_e,1,x_3,\dots,x_n)^{\hat{C}}$ and, consequently, we conclude that
	\begin{align*}
		P_{\mathcal{A}(G)}(x_e,x_3,\dots,x_n)&=\sum_{\hat{C}\in\mathcal{C}(H)}(-1)^{\alpha(\hat{C})}(x_e,1,x_3,\dots,x_n)^{\hat{C}}\\
		&= \sum_{C\in\mathcal{C}(G)}(-1)^{\alpha(C)}(x_e,x_3,\dots,x_n)^C.
	\end{align*}
	If $G$ has multiple loops, by repeating this argument inductively, one can prove the theorem for the most general case.
\end{proof}

We continue by examining the consequences of formula \eqref{eq:charpol} for the spectrum of $D_G$. Unfortunately, a complete characterization of $\sigma(D_G)$ cannot be given. However, it only depends on the structure of $\mathcal{C}(G)$ and the lengths of the cycles present in $G$. The total length of a~subgraph~$H$ of $G$ is given by the sum of the lengths of the edges of $H$, i.e.,
\[L_H=\sum_{e\in E(H)}l_e.\]

\begin{Corollary}
	The characteristic function $P_{\mathcal{A}(G),l}$ of the directed adjacency matrix is given by \looseness=-1
	\[
		P_{\mathcal{A}(G),l}(\lambda)=\sum_{C\in\mathcal{C}(G)}(-1)^{\alpha(C)}\exp({\rm i}\lambda(L_G-L_C)).
	\]
	The spectrum of $D_G$ is given by the solutions to the equation
	\[
		\sum_{C\in\mathcal{C}(G)}(-1)^{\alpha(C)}\exp(-{\rm i}\lambda L_C)=0.
	\]
\end{Corollary}

\begin{proof}
	The first statement follows from a short calculation that uses formula \eqref{eq:charpol}:
	\begin{align*}
		P_{\mathcal{A}(G),l}(\lambda) =\sum_{C\in\mathcal{C}(G)}(-1)^{\alpha(C)}\prod_{e\notin E(C)}\exp({\rm i}\lambda l_e)
		 =\sum_{C\in\mathcal{C}(G)}(-1)^{\alpha(C)}\exp({\rm i}\lambda(L_G-L_C)).
	\end{align*}
	The second statement follows by multiplying the equation $P_{\mathcal{A}(G),l}(\lambda)=0$ by $\exp(-{\rm i}\lambda L_G)$.
\end{proof}

If all edges have length $1$, the characteristic function of the directed adjacency matrix is determined by the characteristic polynomial $p_{\mathcal{A}(G)}$ of $\mathcal{A}(G)$ by
\[
	P_{\mathcal{A}(G),1}(\lambda)= p_{\mathcal{A}(G)}\big({\rm e}^{{\rm i}\lambda}\big)
\]
and we can obtain the characteristic polynomial from the multi-dimensional characteristic polynomial by
\[
	p_{\mathcal{A}(G)}(t)=P_{\mathcal{A}(G)}(t,\dots,t).
\]
This is a well-known object in spectral graph theory (see, e.g., \cite{rigo}), we will prove a few basic facts to illustrate the connection of the characteristic polynomial and hence the spectrum of the Dirac operator to the topology of the graph.

\begin{Corollary}\label{cor:charpol}
	Let $G$ be a directed multigraph with $n$ edges.
	The characteristic polynomial~$p_{\mathcal{A}(G)}$ of the directed adjacency matrix is given by
	\[
		p_{\mathcal{A}(G)}(t)=\sum_{C\in\mathcal{C}(G)}(-1)^{\alpha(C)}t^{n-\eta(C)}.
	\]
\end{Corollary}

\begin{proof}
	The corollary follows directly from formula \eqref{eq:charpol}.
\end{proof}

By expressing the characteristic polynomial in the usual form
\[
	p_{\mathcal{A}(G)}(t)=\sum_{k=0}^n a_k t^k,
\]
Corollary \ref{cor:charpol} implies that the $k$-th coefficient of $p_{\mathcal{A}(G)}(t)$ is given by
\[
	a_k=\sum_{\substack{C\in\mathcal{C}(G) \\ \eta(C)=n-k}}(-1)^{\alpha(C)}.
\]

Recall that the girth $g(G)$ is the size of the shortest directed cycle in $G$ and that the edge-connectivity of $G$ is equal to $k$ if $k$ is the least number of edges that have to be removed such that the remaining graph is disconnected or consists of an isolated vertex. In this case, we say that $G$ is $k$-connected.

\begin{Proposition}
	Let $G$ be a directed multigraph with $n$ edges that contains at least one cycle.
	The girth of $G$ is given by
	\[
		g(G)=\min\{l\in\{1,\dots,n\}\mid a_{n-l}\neq 0\}.
	\]
	Moreover, $-a_{n-l}$ is equal to the number of cycles of size $l$ in $G$ for $0\leq l < 2g(G)$. Furthermore, if $G$ is $k$-connected, then, for $0\leq l <k$, the coefficient $a_l$ is the negative of the number of cycles of size $n-l$ in $G$.
\end{Proposition}

\begin{proof}
	Let $l<g(G)$ and assume that $a_{n-l}$ is not equal to zero. Then, there exists a collection of disjoint cycles $C\in \mathcal{C}(G)$ of size $l$. Consequently, $C$ contains a cycle whose size is less than~$g(G)$, which is a contradiction. Therefore, $a_{n-l}$ must be zero.
	
	If a disjoint collection of cycles $C\in \mathcal{C}(G)$ is of size $l$, where $l$ satisfies $0\leq l < 2g(G)$, then the connectedness of $C$ must be equal to one, otherwise it would contain a cycle of size smaller than $g(G)$. Consequently,
	\[
		a_{n-l}=\sum_{\substack{C\in\mathcal{C}(G) \\ \eta(C)=l}}-1 = -\#\{\text{cycles of size }l\}.
	\]
	Now, let $G$ be $k$-connected and let $0\leq l <k$. If $C\in \mathcal{C}(G)$ is of size $n-l$, then $C$ cannot be disconnected, since only $l$ edges have been removed from $G$ in order to obtain $C$. Therefore, $C$~is a cycle and
	\[
		a_l=\sum_{\substack{C\in\mathcal{C}(G) \\ \eta(C)=n-l}}-1 = -\#\{\text{cycles of size }n-l\}.
\tag*{\qed}
\]
\renewcommand{\qed}{}
\end{proof}

We now characterize some high-order coefficients of the characteristic polynomial.

\begin{Corollary}
	Let $G$ be a graph with $n$ edges.
	For the following coefficients of $p_{\mathcal{A}(G)}$, the following holds:
	\begin{itemize}\itemsep=0pt
		\item $a_n=1$,
		\item $a_{n-1}=-\#\{\text{loops in }G\}$,
		\item $a_{n-2}=-\#\{2\text{-cycles}\}+\#\{\text{pairs of loops that are disjoint}\}$,
		\item $a_{n-3}=-\#\{3\text{-cycles}\} + \#\{\text{pairs of loops and }2\text{-cycles that are disjoint}\}\\ \hphantom{a_{n-3}=} -\#\{\text{triples of loops that are disjoint}\}$.
	\end{itemize}
\end{Corollary}
\begin{proof}
	This follows directly from the above considerations.
\end{proof}

\begin{Example}
	To conclude this section, we give some examples of characteristic polynomials of graphs:
	\begin{itemize}\itemsep=0pt
		\item The characteristic polynomial of the $n$-rose $R_n$ is given by $p_{\mathcal{A}(R_n)}=t^n-nt^{n-1}$.
		\item The characteristic polynomial of the $n$-cycle $C_n$ is given by $p_{\mathcal{A}(C_n)}=t^n-1$.
		\item Let $G_1$ be the graph given in Figure \ref{figure:G1}. The characteristic polynomial of $G_1$ is given by
		$p_{\mathcal{A}(G_1)}=t^4-2t^3$. The reason that the second-order coefficient is equal to zero is due to the fact that the cycle of size $2$ in the middle cancels out the pair of peripheric loops.
		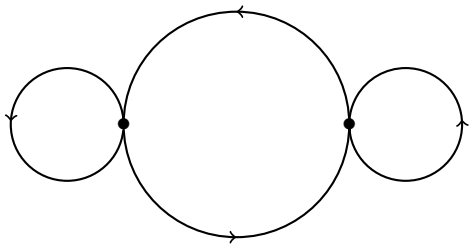
\begin{figure}[H]
			\centering
			\begin{tikzpicture}
				[thick, decoration={
					markings,
					mark=at position 0.5 with {\arrow{>}}}]
				\node (a) at (0,0) {$\bullet$};
				\node (b) at (3,0) {$\bullet$};
				
				\draw[postaction={decorate}] (a)arc(0:354:0.75);
				\draw[postaction={decorate}] (b)arc(180:550:0.75);
				\draw[postaction={decorate}] (a)arc(180:360:1.5);
				\draw[postaction={decorate}] (b)arc(0:180:1.5);
			\end{tikzpicture}
			\caption{The graph $G_1$.}
			\label{figure:G1}
		\end{figure}
		\item Let $G_2$ be the graph given in Figure \ref{figure:G2}. The characteristic polynomial of $G_2$ is given by
		$p_{\mathcal{A}(G_2)}=t^4-2t^3$ and, therefore, is equal to the characteristic polynomial of $G_1$.
		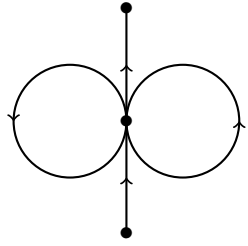
\begin{figure}[H]
			\centering
			\begin{tikzpicture}
				[thick, decoration={
					markings,
					mark=at position 0.5 with {\arrow{>}}}]
				\node (a) at (0,1.5) {$\bullet$};
				\node (b) at (0,0) {$\bullet$};
				\node (c) at (0,-1.5) {$\bullet$};
				
				\draw[postaction={decorate}] (0,-1.5)--(0,0);
				\draw[postaction={decorate}] (0,0)--(0,1.5);
				\draw[postaction={decorate}] (b)arc(0:360:0.75);
				\draw[postaction={decorate}] (b)arc(180:540:0.75);
				
			\end{tikzpicture}
			\caption{The graph $G_2$.}
			\label{figure:G2}
		\end{figure}
		\item Let $G_3$ be the graph given in Figure \ref{figure:G3}. The characteristic polynomial of $G_3$ is given by
		$p_{\mathcal{A}(G_3)}=t^6-3t^4-2t^3$, which, since no two cycles of $G_3$ are disjoint, represents the cycles present in $G_3$, i.e., we have three cycles of size $2$ and two cycles of size $3$.
		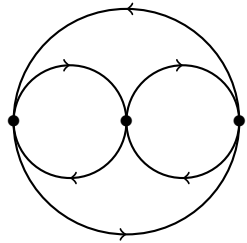
\begin{figure}[H]
			\centering
			\begin{tikzpicture}
				[thick, decoration={
					markings,
					mark=at position 0.5 with {\arrow{>}}}]
				\node (a) at (0,0) {$\bullet$};
				\node (b) at (1.5,0) {$\bullet$};
				\node (c) at (3,0) {$\bullet$};
				
				\draw[postaction={decorate}] (a)arc(180:360:1.5);
				\draw[postaction={decorate}] (c)arc(0:180:1.5);
				\draw[postaction={decorate}] (c)arc(360:180:0.75);
				\draw[postaction={decorate}] (b)arc(360:180:0.75);
				\draw[postaction={decorate}] (a)arc(180:0:0.75);
				\draw[postaction={decorate}] (b)arc(180:0:0.75);
			\end{tikzpicture}
			\caption{The graph $G_3$.}
			\label{figure:G3}
		\end{figure}
	\end{itemize}
\end{Example}

\subsection*{Acknowledgements}

I would like to thank Christian B\"ar, Lashi Bandara and Klaus Ecker for their supervision and for making this project possible. I would also like to thank the anonymous referees for their helpful comments.

\pdfbookmark[1]{References}{ref}
\LastPageEnding


\begin{thebibliography}{99}
\footnotesize\itemsep=0pt

\bibitem{aps}
Atiyah M.F., Patodi V.K., Singer I.M., Spectral asymmetry and {R}iemannian
 geometry.~{I}, \href{https://doi.org/10.1017/S0305004100049410}{\textit{Math. Proc. Cambridge Philos. Soc.}} \textbf{77}
 (1975), 43--69.

\bibitem{bb}
B\"ar C., Ballmann W., Boundary value problems for elliptic differential
 operators of first order, \textit{Surv. Differ. Geom.}, Vol.~17, \href{https://doi.org/10.4310/SDG.2012.v17.n1.a1}{International Press},
 Boston, MA, 2012, 1--78, \href{https://arxiv.org/abs/1101.1196}{arXiv:1101.1196}.

\bibitem{bb2}
B\"ar C., Ballmann W., Guide to elliptic boundary value problems for
 {D}irac-type operators, in Arbeitstagung {B}onn~2013, \textit{Progr. Math.},
 Vol. 319, \href{https://doi.org/10.1007/978-3-319-43648-7_3}{Birkh\"auser}, Cham, 2016, 43--80, \href{https://arxiv.org/abs/1307.3021}{arXiv:1307.3021}.

\bibitem{lashi}
B\"ar C., Bandara L., Boundary value problems for general first-order elliptic
 differential operators, \href{https://doi.org/10.1016/j.jfa.2022.109445}{\textit{J.~Funct. Anal.}} \textbf{282} (2022), 109445,
 69~pages, \href{https://arxiv.org/abs/1906.08581}{arXiv:1906.08581}.

\bibitem{BerKuch}
Berkolaiko G., Kuchment P., Introduction to quantum graphs, \textit{Math. Surv.
 Monogr.}, Vol. 186, \href{https://doi.org/10.1090/surv/186}{American Mathematical Society}, Providence, RI, 2013.

\bibitem{bolte}
Bolte J., Harrison J., Spectral statistics for the {D}irac operator on graphs,
 \href{https://doi.org/10.1088/0305-4470/36/11/307}{\textit{J.~Phys.~A}} \textbf{36} (2003), 2747--2769, \href{https://arxiv.org/abs/nlin.CD/0210029}{arXiv:nlin.CD/0210029}.

\bibitem{bulla}
Bulla W., Trenkler T., The free {D}irac operator on compact and noncompact
 graphs, \href{https://doi.org/10.1063/1.529025}{\textit{J.~Math. Phys.}} \textbf{31} (1990), 1157--1163.

\bibitem{Cvet}
Cvetkovi\'c D.M., Doob M., Sachs H., Spectra of graphs. {T}heory and
 application, VEB Deutscher Verlag der Wissenschaften, Berlin, 1982.

\bibitem{exner}
Exner P., \v{S}eba P., Free quantum motion on a branching graph, \href{https://doi.org/10.1016/0034-4877(89)90023-2}{\textit{Rep.
 Math. Phys.}} \textbf{28} (1989), 7--26.

\bibitem{FullingKuchmentWilson}
Fulling S.A., Kuchment P., Wilson J.H., Index theorems for quantum graphs,
 \href{https://doi.org/10.1088/1751-8113/40/47/009}{\textit{J.~Phys.~A}} \textbf{40} (2007), 14165--14180, \href{https://arxiv.org/abs/0708.3456}{arXiv:0708.3456}.

\bibitem{Harmer}
Harmer M., Hermitian symplectic geometry and extension theory,
 \href{https://doi.org/10.1088/0305-4470/33/50/305}{\textit{J.~Phys.~A}} \textbf{33} (2000), 9193--9203, \href{https://arxiv.org/abs/math-ph/0703027}{arXiv:math-ph/0703027}.

\bibitem{harrisonsurv}
Harrison J.M., Quantum graphs with spin {H}amiltonians, in Analysis on graphs
 and its applications, \textit{Proc. Sympos. Pure Math.}, Vol.~77, \href{https://doi.org/10.1090/pspum/077/2459874}{American
 Mathematical Society}, Providence, RI, 2008, 261--277, \href{https://arxiv.org/abs/0712.0869}{arXiv:0712.0869}.

\bibitem{KostrykinSchrader}
Kostrykin V., Schrader R., Kirchhoff's rule for quantum wires,
 \href{https://doi.org/10.1088/0305-4470/32/4/006}{\textit{J.~Phys.~A}} \textbf{32} (1999), 595--630, \href{https://arxiv.org/abs/math-ph/9806013}{arXiv:math-ph/9806013}.

\bibitem{Kuchment2014}
Kuchment P., Quantum graphs:~{I}.~{S}ome basic structures, \textit{Waves Random
 Media} \textbf{14} (2004), S107--S128.

\bibitem{kuchment2008quantum}
Kuchment P., Quantum graphs: an introduction and a brief survey, in Analysis on
 {G}raphs and its {A}pplications, \textit{Proc. Sympos. Pure Math.}, Vol.~77,
 \href{https://doi.org/10.1090/pspum/077/2459876}{American Mathematical Society}, Providence, RI, 2008, 291--312,
 \href{https://arxiv.org/abs/0802.3442}{arXiv:0802.3442}.

\bibitem{post}
Post O., First order approach and index theorems for discrete and metric
 graphs, \href{https://doi.org/10.1007/s00023-009-0001-3}{\textit{Ann. Henri Poincar\'e}} \textbf{10} (2009), 823--866,
 \href{https://arxiv.org/abs/0708.3707}{arXiv:0708.3707}.

\bibitem{rigo}
Rigo M., Advanced graph theory and combinatorics, \textit{Comput. Eng. Ser.}, \href{https://doi.org/10.1002/9781119008989}{ISTE},
 London, 2016.

\bibitem{VONBELOW}
von Below J., A characteristic equation associated to an eigenvalue problem on
 {$c^2$}-networks, \href{https://doi.org/10.1016/0024-3795(85)90258-7}{\textit{Linear Algebra Appl.}} \textbf{71} (1985), 309--325.

\end{thebibliography}
\end{document}